\documentclass[twocolumn]{autart}    

\usepackage{graphicx}          

\usepackage{mathtools}
\usepackage{amsmath}
\usepackage{units}

\usepackage{cite}
\usepackage{amssymb} 
\usepackage{color}
\newcommand{\T}{\scriptscriptstyle\top}  
\newcommand{\defeq}{\vcentcolon=}
\usepackage{enumerate}
\usepackage{enumitem}
\DeclareUnicodeCharacter{2212}{\ensuremath{-}}

\makeatletter
\let\save@mathaccent\mathaccent
\newcommand*\if@single[3]{%
  \setbox0\hbox{${\mathaccent"0362{#1}}^H$}%
  \setbox2\hbox{${\mathaccent"0362{\kern0pt#1}}^H$}%
  \ifdim\ht0=\ht2 #3\else #2\fi
  }
\newcommand*\rel@kern[1]{\kern#1\dimexpr\macc@kerna}
\newcommand*\widebar[1]{\@ifnextchar^{{\wide@bar{#1}{0}}}{\wide@bar{#1}{1}}}
\newcommand*\wide@bar[2]{\if@single{#1}{\wide@bar@{#1}{#2}{1}}{\wide@bar@{#1}{#2}{2}}}
\newcommand*\wide@bar@[3]{%
  \begingroup
  \def\mathaccent##1##2{%
    \let\mathaccent\save@mathaccent
    \if#32 \let\macc@nucleus\first@char \fi
    \setbox\z@\hbox{$\macc@style{\macc@nucleus}_{}$}%
    \setbox\tw@\hbox{$\macc@style{\macc@nucleus}{}_{}$}%
    \dimen@\wd\tw@
    \advance\dimen@-\wd\z@
    \divide\dimen@ 3
    \@tempdima\wd\tw@
    \advance\@tempdima-\scriptspace
    \divide\@tempdima 10
    \advance\dimen@-\@tempdima
    \ifdim\dimen@>\z@ \dimen@0pt\fi
    \rel@kern{0.6}\kern-\dimen@
    \if#31
      \overline{\rel@kern{-0.6}\kern\dimen@\macc@nucleus\rel@kern{0.4}\kern\dimen@}%
      \advance\dimen@0.4\dimexpr\macc@kerna
      \let\final@kern#2%
      \ifdim\dimen@<\z@ \let\final@kern1\fi
      \if\final@kern1 \kern-\dimen@\fi
    \else
      \overline{\rel@kern{-0.6}\kern\dimen@#1}%
    \fi
  }%
  \macc@depth\@ne
  \let\math@bgroup\@empty \let\math@egroup\macc@set@skewchar
  \mathsurround\z@ \frozen@everymath{\mathgroup\macc@group\relax}%
  \macc@set@skewchar\relax
  \let\mathaccentV\macc@nested@a
  \if#31
    \macc@nested@a\relax111{#1}%
  \else
    \def\gobble@till@marker##1\endmarker{}%
    \futurelet\first@char\gobble@till@marker#1\endmarker
    \ifcat\noexpand\first@char A\else
      \def\first@char{}%
    \fi
    \macc@nested@a\relax111{\first@char}%
  \fi
  \endgroup
}

\usepackage{float}

\usepackage{url}
\allowdisplaybreaks
\usepackage{hyperref}

\begin{document}
\date{}
\begin{frontmatter}

\title{Synthesis of Discrete-time Control Barrier Functions for Polynomial Systems Based on Sum-of-Squares Programming\thanksref{footnoteinfo}} 

\thanks[footnoteinfo]{The research leading to these results has received funding from the European Research Council under the Advanced ERC Grant Agreement PROACTHIS, no. 101055384.}

\author[Eindhoven]{Erfan Shakhesi}\ead{e.shakhesi@tue.nl},    
\author[Eindhoven]{W.P.M.H. (Maurice) Heemels}\ead{m.heemels@tue.nl},               
\author[Eindhoven]{Alexander Katriniok}\ead{a.katriniok@tue.nl}  

\address[Eindhoven]{Department of Mechanical Engineering, Control Systems Technology Section, Eindhoven
University of Technology, \mbox{The Netherlands}.}  


\begin{abstract}
Discrete-time Control Barrier Functions (DTCBFs) are commonly utilized in the literature as a powerful tool for synthesizing control policies that guarantee safety of discrete-time dynamical systems. However, the systematic synthesis of DTCBFs in a computationally efficient way is at present an important open problem. This article first proposes a novel alternating-descent approach based on Sum-of-Squares programming to synthesize quadratic DTCBFs and corresponding polynomial control policies for discrete-time control-affine polynomial systems with input constraints and semi-algebraic safe sets. Subsequently, two distinct approaches are introduced to extend the proposed method to the synthesis of higher-degree polynomial DTCBFs. To demonstrate its efficacy, we apply the proposed method to numerical case studies.
\end{abstract}

\end{frontmatter}

\section{Introduction} \label{sec:sec1}
Designing controllers that guarantee safety is essential in many applications, such as autonomous driving, where safety violations can lead to severe consequences. In the context of dynamical systems, a system is often considered safe as long as its trajectories remain within a predefined safe set. Under this definition, safety can be guaranteed if all initial states lie within a controlled invariant set that is itself a subset of the safe set. To construct such controlled invariant sets and synthesize safe controllers, Control Barrier Functions (CBFs) \cite{Ames2014a} have become widely used in the literature.

CBFs are Lyapunov-like functions that guarantee the controlled invariance of their zero-superlevel set, defined as the set of all states for which a CBF is non-negative. Using such a function in an online optimization problem, safe controllers can be synthesized by adjusting known nominal controllers when they lead to a violation of safety constraints; in this sense, CBFs can be used to act as a safety filter. This CBF-based control technique was initially proposed for continuous-time systems \cite{Ames2014a} and has intimate connections to projected dynamical systems \cite{Giannis-TAC}.
An extension of CBFs to discrete-time dynamical systems, referred to as discrete-time CBFs or DTCBFs for short, has been introduced in \cite{Agrawal2017a}. DTCBFs have potential for various applications, including their use within Model Predictive Control frameworks \cite{katriniok2023, Zeng2021a}. 

The synthesis of CBFs for continuous-time systems has gained considerable attention recently, with methods such as Sum-of-Squares (SOS) programming \cite{Ames2019a, dai2022, clark2022, Wang2023b, Zhang2023, Kang2023, wang2024convexcodesign}, backup policy methods \cite{Gurriet2020, chen2021}, and sampling-based approaches \cite{wei2023, Zhang2023, lindemann2024}.
However, the systematic synthesis of DTCBFs is more challenging, and research in this area remains limited. Notably, prior works that exploited DTCBFs, such as \cite{Agrawal2017a}, \cite{Zeng2021b}, and \cite{Zeng2021a}, have resorted to utilizing \textit{candidate} DTCBFs. These candidates, however, may not rigorously adhere to the formal DTCBF definition, particularly in the presence of input constraints, and therefore do not guarantee safety and feasibility. Nevertheless, \cite{Zhang2024} proposes a learning-based method for synthesizing candidate DTCBFs by minimizing safety violations on sample states. Additionally, \cite{Freire2023a, RDTCBF-literature} leverage maximal output admissible set (MOAS) theory to synthesize a particular type of DTCBF (considering the class-$\mathcal{K}_{\infty}$ function as the identity in the formal DTCBF definition). The MOAS theory, though, has not been well-studied for nonlinear systems. Moreover, SOS programming has been used in \cite{Wang2023a, Zamani-2020} to synthesize polynomial DTCBFs for polynomial systems, though with some limitations that we review below.

The conditions required to synthesize polynomial (DT)CBFs for polynomial systems and semi-algebraic safe sets can be mapped to SOS constraints via the Positivstellensatz or the generalized S-procedure lemma \cite{Ames2019a}. Due to bilinearity, these constraints cannot be imposed directly in SOS programming. Existing research on \textit{continuous-time} systems has addressed this issue for control-affine systems by alternating-descent approaches \cite{Ames2019a, dai2022, clark2022, Wang2023a, Zhang2023}. 

Unfortunately, the DTCBF constraint is generally not affine in the control input, even if the system is control-affine. This introduces additional challenges in synthesizing DTCBFs using SOS programming, compared to the continuous-time case, due to the presence of bilinear terms in the DTCBF constraint. 
As illustrated in \mbox{Example \ref{example:sec4} below}, these bilinear terms cannot be addressed using the alternating-descent approaches proposed in \cite{Ames2019a, dai2022, clark2022, Wang2023a, Zhang2023}. The same issue exists for synthesizing discrete-time Control Lyapunov Functions (CLFs) and \textit{classical} controlled invariant sets for discrete-time systems, except for the synthesis of polyhedral controlled invariant sets considering control-affine systems. To our knowledge, addressing these bilinear terms using SOS programming and, in turn, synthesizing DTCBFs, discrete-time CLFs, and controlled invariant sets\textemdash aside from polyhedral sets\textemdash for discrete-time polynomial systems remains an open research problem in the literature.

Nonetheless, \cite{Wang2023a, Zamani-2020} have attempted to bypass the issue of bilinearity to synthesize DTCBFs. Specifically, \cite{Wang2023a} proposes an alternating-descent approach that can be exploited to synthesize polyhedral controlled invariant sets or (multiple) \textit{affine} DTCBFs for control-affine systems. Although this approach is interesting, affine DTCBFs have limitations in applicability, and polyhedral controlled invariant sets are not suitable for applications such as obstacle avoidance when resulting controlled invariant sets must be non-convex. Alternatively, in \cite{Zamani-2020}, nonlinear polynomial DTCBFs are synthesized using a set-valued control policy rather than a specific polynomial as a corresponding control policy. In particular, this approach assumes that all control inputs that are element-wise greater than or equal to a vector of polynomials form the corresponding control policy.
This is noteworthy, but it may not be suitable for certain practical applications, especially when the control input at specific states is required to take only negative values. 

In our preliminary conference paper \cite{Shakhesi2025}, we addressed the aforementioned challenges associated with quadratic DTCBFs, but without the technical proofs due to space limitations. In this article, we provide complete proofs and also introduce two distinct approaches to extend the proposed synthesis method to higher-degree polynomial DTCBFs and corresponding polynomial control policies for control-affine polynomial systems with input constraints and semi-algebraic safe sets. Our proposed method is also of interest in the context of synthesizing controlled invariant sets and CLFs for discrete-time polynomial systems. Furthermore, we apply the proposed method to various numerical case studies, including a novel system with a higher dimension than the one considered in \cite{Shakhesi2025}.

\textit{Notation:} We use $\mathbb{R}$, $\mathbb{R}_{> 0}$, and $\mathbb{R}_{\geqslant 0}$ to denote the set of real numbers, positive real numbers, and non-negative real numbers, respectively. Additionally, $\mathbb{R}^{n}$ and $\mathbb{R}^{n \times n}$ represent the set of all $n$-dimensional vectors and \mbox{$n \times n$} matrices, respectively, of real numbers. Similarly, $\mathbb{R}^{n}_{\geqslant 0}$ represents the set of all $n$-dimensional vectors of non-negative real numbers.
For a vector \mbox{$x \in \mathbb{R}^n$}, $x_i \in \mathbb{R}$, \mbox{$i \in \{1,2, \hdots, n\}$}, represents the $i$-th element of $x$, and $x^{\T}$ represents its transpose. Moreover, $\mathbb{N} \defeq \{1,2,3,\hdots\}$, $\mathbb{N}_0 \defeq \mathbb{N} \cup \{0\}$. 
For a set $\mathcal{S} \subset \mathbb{R}^n$, we use $\mathcal{S}^{\mathrm{c}}$ and $\partial \mathcal{S}$ to represent its complement and its boundary, respectively.
The zero-superlevel set $\mathcal{C}$ of a function $h:\mathbb{R}^{n} \rightarrow \mathbb{R}$ is defined as 
    $\mathcal{C} \defeq \{ x \in \mathbb{R}^{n} \mid h(x) \geqslant 0 \}$. 
    A continuous function \mbox{$\gamma:\mathbb{R}_{\geqslant 0} \rightarrow \mathbb{R}_{\geqslant 0}$} is said to be a class-$\mathcal{K}_\infty$ function, denoted by $\gamma \in \mathcal{K}_{\infty}$, if it is strictly increasing, $\gamma(0) = 0$, and $\lim_{r\rightarrow \infty} \gamma(r)=\infty$. For $\gamma \in \mathcal{K}_{\infty}$, we use the notation $\gamma \in \mathcal{K}^{\leqslant \mathrm{id}}_{\infty}$ to indicate that $\gamma(r) \leqslant r$ for all $r \in \mathbb{R}_{\geqslant 0}$.
   
\section{Background on Discrete-time Control Barrier Functions} \label{sec:sec2}
\vspace{-0.2cm}
Consider the control-affine discrete-time system 
\begin{equation} \label{eq:sec2:dynamical-system}
    x^+ = f(x)+g(x)u, 
\end{equation}
with state vector $x \in \mathbb{R}^n$, control input vector \mbox{$u \in \mathbb{U} \subseteq \mathbb{R}^m$}, both at the current time instant, state vector $x^+ \in \mathbb{R}^n$ at the next time instant, and the mappings $f:\mathbb{R}^n \rightarrow \mathbb{R}^n$ and $g:\mathbb{R}^n \rightarrow \mathbb{R}^{n \times m}$. Here, $\mathbb{U}$ is the control admissible set. Moreover, the safe set $\mathcal{S}$ for the system \eqref{eq:sec2:dynamical-system} is defined as
\begin{align} \label{eq:sec2:safe-set}
    \mathcal{S} \defeq \{ x \in \mathbb{R}^n \mid s(x) \geqslant 0\},
\end{align}
where $s:\mathbb{R}^n \rightarrow \mathbb{R}$ is a given mapping.
\begin{defn}[Controlled invariance]
For the system \eqref{eq:sec2:dynamical-system} with the control admissible set $\mathbb{U}$, a set $\mathcal{C} \subset \mathbb{R}^n$ is controlled invariant, if for all $x \in \mathcal{C}$, there exists a control input $u \in \mathbb{U}$ such that $f(x)+g(x)u \in \mathcal{C}$.
\end{defn}
\begin{defn}[Safety] \label{def:sec2:safety}
    The system \eqref{eq:sec2:dynamical-system} with the control admissible set $\mathbb{U}$ is considered safe with respect to a safe set $\mathcal{S} \subset \mathbb{R}^n$ and an initial state set \mbox{$X_0 \subset \mathbb{R}^n$}, denoted by $(\mathcal{S}, X_0)$-safe, if there exists a controlled invariant \mbox{set $\mathcal{C} \subset \mathbb{R}^n$} such that $X_0 \subseteq \mathcal{C} \subseteq \mathcal{S}$.
\end{defn}
To guarantee safety under Definition \ref{def:sec2:safety}, DTCBFs are commonly utilized in the literature, which also enables the synthesis of safe controllers. 
\begin{defn}[DTCBF \cite{Agrawal2017a}, \cite{Zeng2021a}] \label{def:sec2:DTCBF}
{\sloppy
     Consider a function \mbox{$h:\mathbb{R}^n \rightarrow \mathbb{R}$} with zero-superlevel set $\mathcal{C}$. For the \mbox{system \eqref{eq:sec2:dynamical-system}} with the control admissible set $\mathbb{U}$, $h$ is a \textit{discrete-time Control Barrier Function (DTCBF)}, if there exists a \mbox{$\gamma \in \mathcal{K}^{\leqslant \mathrm{id}}_{\infty}$} such that for each $x \in \mathcal{C}$, there exists a control input $u \in \mathbb{U}$ satisfying
    \begin{align} \label{eq:sec2:DTCBF-constraint-alpha}
        h(f(x)+g(x)u) - h(x) \geqslant - \gamma(h(x)).
    \end{align}
    Additionally, we say that a control policy $\pi : \mathcal{C} \rightarrow \mathbb{R}^{m}$ is a \textit{friend} of a DTCBF $h$ and a corresponding \mbox{$\gamma \in \mathcal{K}^{\leqslant \mathrm{id}}_{\infty}$}, denoted by the DTCBF-triple $(h, \gamma, \pi)$, for the system \eqref{eq:sec2:dynamical-system} with $\mathbb{U}$, if for all $x \in \mathcal{C}$,
    \begin{align}
        \pi(x) \in \mathbb{U} \text{~~and~~}
        h(f(x)+ g(x) \pi(x)) - h(x) \geqslant - \gamma(h(x)). \nonumber
    \end{align}
}
\end{defn}

\begin{thm}[Safety \cite{Agrawal2017a}] \label{theorem:sec2:safety}
     For the system \eqref{eq:sec2:dynamical-system} with the control admissible set \hspace{0.05cm}$\mathbb{U}$, consider a DTCBF $h$ with zero-superlevel set $\mathcal{C}$. Then, $\mathcal{C}$ is \textit{controlled invariant}, and the \mbox{system \eqref{eq:sec2:dynamical-system}} with $\mathbb{U}$ is \textit{$(\mathcal{S}, X_0)$-safe}, if \mbox{$X_0 \subseteq \mathcal{C} \subseteq \mathcal{S}$}.
\end{thm}

\begin{rem}
    In addition to the fact that the zero-superlevel set $\mathcal{C}$ of a DTCBF $h$ is \textit{controlled invariant}, the DTCBF constraint \eqref{eq:sec2:DTCBF-constraint-alpha} includes an additional term, $\gamma \in \mathcal{K}^{\leqslant \mathrm{id}}_{\infty}$. This term regulates the rate at which the states of the system \eqref{eq:sec2:dynamical-system} can approach the \mbox{boundary of $\mathcal{C}$}. The advantages of adjusting $\gamma$ and using DTCBFs in Model Predictive Control frameworks are discussed in more detail in \cite{katriniok2023} and \cite{Zeng2021a}. If $\gamma$ is chosen to be the identity function, the DTCBF definition is equivalent to the definition of controlled invariance. 
\end{rem}


\section{Background on SOS Programming} \label{sec:sec3}
In this section, we provide a brief overview of SOS programming. For more details, we refer the reader to \cite{sostools} and \cite{SOS2005}. 
\begin{defn}[SOS polynomial]
    A polynomial $p$ in $x \in \mathbb{R}^n$ with real coefficients, denoted by \mbox{$p \in \mathbb{R}[x]$}, is called a Sum-of-Squares (SOS) polynomial, denoted by \mbox{$p \in \mathrm{SOS}[x]$}, if there exist $q_i \in \mathbb{R}[x]$, \mbox{$i \in \{1,2, \hdots, I\}$}, $I \in \mathbb{N}$, such that
        \mbox{$p(x) = \sum_{i = 1}^{I} q^2_i(x)$} for all $x \in \mathbb{R}^n$. 
\end{defn}
Moreover, we denote $\mathbb{R}[x]^r$ and $\mathbb{R}[x]^{r \times r}$ as the sets of $r$-dimensional polynomials $p: \mathbb{R}^n \rightarrow \mathbb{R}^r$ and $r \times r$-dimensional polynomials $q: \mathbb{R}^n \rightarrow \mathbb{R}^{r \times r}$, respectively, in \mbox{$x \in \mathbb{R}^n$} with real coefficients. Similarly, we denote $\mathrm{SOS}[x]^r$ as the set of $r$-dimensional SOS polynomials.

For a polynomial $p \in \mathbb{R}[x]$, the condition that $p$ is a SOS polynomial (\textit{SOS constraint}), $p \in \mathrm{SOS}[x]$, implies that
    $p(x) \geqslant 0$ for all $x \in \mathbb{R}^n$,
but the converse does not hold generally \cite{Hilbert1888}. Nonetheless, verifying a SOS constraint is computationally more efficient than verifying the non-negativity of a polynomial. 
Non-negativity conditions within semi-algebraic sets can be guaranteed by SOS constraints using the generalized S-procedure lemma below.
\begin{lem}[Generalized S-procedure \cite{S-procedure, Wang2023a}] \label{lemma:sec4:S-procedure}
    Consider a semi-algebraic set
    \begin{align}
        \mathcal{W} \defeq \{x \in \mathbb{R}^n &\mid w_i(x) \geqslant 0, ~ i \in \{1,2, \hdots, I\}\}, \nonumber
    \end{align} 
    where $w_i\in \mathbb{R}[x]$, $i \in \{1,2, \hdots, I\}$, and $I \in \mathbb{N}$. The condition \mbox{$l(x) \geqslant 0$} for an $l \in \mathbb{R}[x]$ holds for all $x \in \mathcal{W}$, if there exist \mbox{$p_i \in \mathrm{SOS}[x]$}, $i \in \{1,2, \hdots, I\}$, such that \vspace{-0.2cm}
    \begin{align} 
        L(x) \defeq l(x) - \sum_{i=1}^I p_i(x)w_i(x) \geqslant 0 \text{~~ for all $x \in \mathbb{R}^n$} \nonumber
    \end{align} \\ [-0.2cm]
    or if $L \in \mathrm{SOS}[x]$.
    Additionally, $l(x) > 0$ holds for all \mbox{$x \in \mathcal{W}$}, if there exist $\epsilon \in \mathbb{R}_{>0}$ and \mbox{$p_i \in \mathrm{SOS}[x]$}, \mbox{$i \in \{1,2, \hdots, I\}$}, such that \vspace{-0.2cm} 
    \begin{align}
        L'(x) \defeq l(x) - \epsilon - \sum_{i=1}^I p_i(x)w_i(x) \geqslant 0 \text{~~ for all $x \in \mathbb{R}^n$} \nonumber
    \end{align} \\ [-0.2cm]
    or if $L' \in \mathrm{SOS}[x]$.
\end{lem}

A \textit{standard} feasibility problem in SOS programming involves finding real-valued decision variables that satisfy SOS constraints, where the decision variables appear linearly. This kind of problem can be efficiently solved using semidefinite programming \cite{sostools}. 

\section{Problem Statement} \label{sec:sec4}
\vspace{-0.2cm}
In this section, we discuss the synthesis problem, the assumptions considered to address it, and the main challenge of solving this problem using SOS programming.
\begin{prob}[Synthesis of DTCBFs] \label{prob:sec3:synthesis}
    Consider the system \eqref{eq:sec2:dynamical-system} with the control admissible set $\mathbb{U}$ and the safe set $\mathcal{S}$ in \eqref{eq:sec2:safe-set}. The objective is to synthesize a DTCBF-triple $(h,\gamma, \pi)$ according to Definition \ref{def:sec2:DTCBF} such that the size of the zero-superlevel set $\mathcal{C}$ of $h$ is maximized (in an appropriate sense), and $\mathcal{C} \subseteq \mathcal{S}$.
\end{prob}
\begin{assum} \label{assump:sec3:dynamics-polynomial}
    The mappings $f$ and $g$ associated with the system \eqref{eq:sec2:dynamical-system} and the mapping $s$ associated with the safe set $\mathcal{S}$ in \eqref{eq:sec2:safe-set} are polynomials with real coefficients. 
\end{assum}
\vspace{-0.2cm}
\begin{assum} \label{assump:sec3:m-rec-control-admissible}
    The control admissible set $\mathbb{U}$ is polyhedral in the sense that it can be \mbox{defined as} \vspace{-0.2cm}
    \begin{align} \label{eq4:sec4:control-admissible-assumption}
        \mathbb{U} \defeq \{ u \in \mathbb{R}^m \mid Mu + d \geqslant 0\},
    \end{align} \\[-0.5cm]
    for some $M \in \mathbb{R}^{n_U \times m}$ and $d \in \mathbb{R}^{n_U}$. 
    
\end{assum}
To address Problem \ref{prob:sec3:synthesis} under Assumptions \ref{assump:sec3:dynamics-polynomial} and \ref{assump:sec3:m-rec-control-admissible} using SOS programming, we choose to search for a polynomial DTCBF $h$ with polynomials $\pi$ and $\gamma$. In particular, we parameterize $h$ as
\begin{align} \label{eq:sec4:parameterized-DTCBF-general}
    h(x) \defeq \vartheta^{\T}H(x), 
\end{align}
where $H \in \mathbb{R}[x]^{n_{\vartheta}}$ is a user-specified vector of monomials, and $\vartheta \in \mathbb{R}^{n_{\vartheta}}$ is the vector of unknown coefficients.
Additionally, we parameterize $\gamma \in \mathcal{K}^{\leqslant \mathrm{id}}_{\infty}$ as a linear polynomial, which is commonly utilized in the literature \cite{Agrawal2017a, Zeng2021a, katriniok2023, shakhesi2024b},
\begin{align} \label{eq:sec4:parameterized-gamma}
    \gamma (r) \defeq \gamma_0r, \quad r\in \mathbb{R}_{\geqslant 0}, 
\end{align}
with parameter $\gamma_0 \in (0,1] \subset \mathbb{R}$ unknown.
Likewise, we parameterize the control policy $\pi$ as
\begin{align} \label{eq:sec4:parameterized-control-policy}
    \pi_i(x) \defeq \kappa_{i}^{\T}\Pi_{i}(x), 
\end{align}
$i \in \{1,2, \hdots, m\}$, where $\pi(x) \defeq [\pi_1(x)~ \hdots~ \pi_m(x)]^{\T}$. Here, $\Pi_{i} \in \mathbb{R}[x]^{n_{\kappa_i}}$, $\kappa_{i} \in \mathbb{R}^{n_{\kappa_i}}$, and $n_{\kappa_i} \in \mathbb{N}$ represent a user-specified vector of monomials, the vector of unknown coefficients, and the number of coefficients, respectively, in $\pi_i$, $i \in \{1,2, \hdots, m\}$.

Based on Definition \ref{def:sec2:DTCBF}, $(h, \gamma, \pi)$, parameterized as in \mbox{\eqref{eq:sec4:parameterized-DTCBF-general}--\eqref{eq:sec4:parameterized-control-policy}}, is a DTCBF-triple, if and only if for all $x \in \mathcal{C}$, where $\mathcal{C}$ is the zero-superlevel set of $h$, it holds that
\begin{align}
    h(f(x) + g(x) \pi(x)) -  
    h(x) + \gamma(h(x)) & \geqslant 0, \label{eq:sec5:conditions-DTCBF-1}\\
    M\pi(x) + d & \geqslant 0. \label{eq:sec5:conditions-DTCBF-2}
\end{align}
Additionally, we aim to ensure that the condition $\mathcal{C} \subseteq \mathcal{S}$ holds (see Theorem \ref{theorem:sec2:safety}), thereby making the system \eqref{eq:sec2:dynamical-system} with $\mathbb{U}$ $(\mathcal{S}, X_0)$-safe, where $X_0 \subseteq \mathcal{C}$. This condition is satisfied, if and only if it holds that $s(x) \hspace{-0.03cm} \geqslant \hspace{-0.03cm} 0$ for all $x \in \mathcal{C}$, where $s$ is associated with the safe set $\mathcal{S}$ as \mbox{in \eqref{eq:sec2:safe-set}}. However, given that $s$ is known and $h$ is unknown, expressing this condition as a SOS constraint using the generalized S-procedure (Lemma \ref{lemma:sec4:S-procedure}) leads to a bilinearity issue that can be avoided. Specifically, we can reformulate it by stating that $\mathcal{C} \subseteq \mathcal{S}$ holds, if and only if 
\begin{align}
    h(x) &< 0 \qquad \text{for all } x \in \mathcal{S}^{\mathrm{c}} \defeq \mathbb{R}^n \backslash \mathcal{S}. \label{eq:sec5:conditions-DTCBF-3}
\end{align}
Similar to \cite{Wang2023a, Ames2019a} and based on the generalized S-procedure, 
we can state that the conditions \eqref{eq:sec5:conditions-DTCBF-1}\textendash \eqref{eq:sec5:conditions-DTCBF-3} hold, if there exist $\Omega, \Phi \in \text{SOS}[x]$, \mbox{$\Psi \in \text{SOS}[x]^{n_U}$}, and $\epsilon \in \mathbb{R}_{>0}$ such that
\begin{subequations}
\begin{align}
    h(f(x) + g(x) \pi(x)) − h(x) + \nonumber \\
    \gamma(h(x)) −\Omega(x)h(x) &\in \text{SOS}[x], \label{eq:sec4:required-SOS-2}\\
    M\pi(x) + d − \Psi(x)h(x) &\in \text{SOS}[x]^{n_U}, \label{eq:sec4:required-SOS-3} \\
    -h(x) - \epsilon + \Phi(x)s(x) &\in \text{SOS}[x]. \label{eq:sec4:required-SOS-1}
\end{align}
\end{subequations}
Unfortunately, the presence of bilinear terms (products of decision variables), such as $\gamma(h(x))$, $\Omega(x)h(x)$, $\Psi(x)h(x)$, and $h(f(x)+g(x)\pi(x))$, poses a challenge as they cannot be directly imposed in a standard feasibility problem in SOS programming, where decision variables appear linearly. For synthesizing \textit{continuous-time} CBFs, alternating-descent approaches are proposed in \cite{clark2022, Wang2023a, Ames2019a, dai2022} to handle bilinear terms that resemble $\gamma(h(x))$, $\Omega(x)h(x)$, and $\Psi(x)h(x)$.
Specifically, in these approaches, one of the two polynomials in a bilinear term is kept fixed, while the other is updated. Subsequently, the one that was obtained becomes fixed, and the previously fixed polynomial is updated. However, the bilinearity of $h(f(x)+g(x)\pi(x))$ cannot be addressed using these approaches. To make this clear, consider the following motivating example.
\begin{exmp}[Motivating example] \label{example:sec4}
Consider a simple dynamical system
\begin{align}
    x^+ = x + u, \nonumber
\end{align}
where $x \in \mathbb{R}$ and $u \in \mathbb{U} \subseteq \mathbb{R}$. Assume that we parameterize a DTCBF $h$ and a control policy $\pi$ as 
\begin{align}
    h(x) &\defeq \vartheta_0 + \vartheta_1 x + \vartheta_2 x^2, \nonumber \\
    \pi(x) &\defeq \kappa_{1} + \kappa_{2} x, \nonumber
\end{align}
where $\vartheta_0, \vartheta_1, \vartheta_2, \kappa_{1}, \kappa_{2} \in \mathbb{R}$ are unknown. The value of the DTCBF at the next time instant is 
\vspace{-0.18cm}
\begin{align} \label{eq:motivating-example:DTCBF-next}
    h(x + \pi(x)) &= \vartheta_0 + \vartheta_1 \bigl( x + \kappa_1 + \kappa_2 x \bigr) + \vartheta_2 [ (1 + 2\kappa_{2} + \nonumber \\
    &  \quad ~~~ \kappa_{2}^2)x^2 + (2\kappa_{1} + 2\kappa_{1} \kappa_{2})x + \kappa_{1}^2].
\end{align} \\[-0.5cm]
As a result, in the SOS constraint \eqref{eq:sec4:required-SOS-2},  we encounter bilinear terms involving $\kappa_{2}^2$, $\kappa_{1} \kappa_{2}$, and $\kappa_{1}^2$ even by applying alternating-descent approaches\textemdash where we fix $h$ (in this case $\vartheta_0$, $\vartheta_1$, and $\vartheta_2$) and aim to obtain the unknown coefficients of $\pi$ (in this case $\kappa_{1}$ and $\kappa_{2}$). Note that by choosing $\vartheta_2 \defeq 0$, there is no bilinearity issue in \eqref{eq:motivating-example:DTCBF-next} when using a standard alternating-descent approach, which allows the synthesis of affine DTCBFs and polyhedral controlled invariant sets, as done in \cite{Wang2023a}.
\end{exmp}

\begin{rem}
    Note that the zero-superlevel set of affine DTCBFs is polyhedral. As a result, it may not be suitable for certain nonlinear systems in constructing controlled invariant sets or guaranteeing safety, particularly in cases where the controlled invariant set must be non-convex, such as in obstacle avoidance applications. 
\end{rem}
Our main contribution in this article is to propose a method to handle the bilinear terms in $h(f(x) + g(x) \pi(x))$. 
To simplify the introduction of the proposed method, we first explain it in Section \ref{sec:sec5} for quadratic DTCBFs $h$ parameterized as 
\begin{align} \label{eq:sec4:parameterized-DTCBF}
    h(x) \defeq \vartheta_0 + \vartheta_1^{\T}x + x^{\T}\vartheta_2x, 
\end{align} 
where $\vartheta_0 \in \mathbb{R}$, $\vartheta_1 \in \mathbb{R}^n$, and \mbox{$\vartheta_2 \in \mathbb{R}^{n \times n}$} is a symmetric matrix, all of which are unknown. 
Then, in Section \ref{sec:sec6}, we extend the proposed method to synthesize polynomial DTCBFs of higher degrees. 
\section{Proposed Synthesis Method for Quadratic DTCBFs Based on SOS Programming} \label{sec:sec5}
The proposed alternating-descent approach for quadratic DTCBFs consists of three steps: initialization (\mbox{Step \ref{step:step0}}), updating the control policy (Step \ref{step:step1}), and updating \mbox{the DTCBF (Step \ref{step:step2})}. To be more precise, in Step \ref{step:step0}, we initialize \mbox{the DTCBF $h^{(0)}$}. In particular, $h^{(0)}$ needs to be a valid quadratic DTCBF but it can be quite conservative. Using $h^{(0)}$ in \mbox{Step \ref{step:step1}}, we synthesize $\gamma^{(1)}$ and its \mbox{friend $\pi^{(1)}$}. Subsequently, with $\gamma^{(1)}$ and $\pi^{(1)}$ fixed, we obtain the updated DTCBF $h^{(1)}$ in \mbox{Step \ref{step:step2}} such that $\mathcal{C}^{0} \subset \mathcal{C}^{(1)}$, where $\mathcal{C}^{(0)}$ and $\mathcal{C}^{(1)}$ are the zero-superlevel sets of \mbox{$h^{(0)}$ and $h^{(1)}$}, respectively. We then proceed to Step \ref{step:step1}. This procedure continues until the size of the zero-superlevel set can no longer be enlarged. The main challenge in this procedure occurs in \mbox{Step \ref{step:step1}} due to the bilinearity issue (see \mbox{Example \ref{example:sec4}).} In \mbox{Section \ref{sec:sec5.1}}, we first propose a method to handle this issue for quadratic DTCBFs, and then we explain the algorithm overview in Section \ref{sec:sec5.2}.

\vspace{-0.1cm}
\subsection{Handling the Bilinearity Issue in Step \ref{step:step1}} \label{sec:sec5.1}
\vspace{-0.1cm}
Given the quadratic DTCBF $h^{(k-1)}$ obtained at the $(k-1)$-th iteration, $k \in \mathbb{N}$, to update $\gamma$ and its friend $\pi$ for the system \eqref{eq:sec2:dynamical-system} with the control admissible set $\mathbb{U}$ at the $k$-iteration, we need to impose \vspace{-0.3cm} 
\begin{subequations} 
\begin{align}
    &h^{(k-1)}\bigl(\overbrace{f(x) + g(x)\pi(x)}^{x^+ \defeq}\bigr) - h^{(k-1)}\left(x\right) + \nonumber \\
    &\hspace{0.645cm} \gamma(h^{(k-1)}(x)) - \Omega(x) h^{(k-1)}(x) \in \text{SOS}[x], \label{eq:sec5:updating_CBF_original_function} \\
    &~~~~~~~~M\pi(x) + d -\Psi(x) h^{(k-1)}(x) \in \text{SOS}[x]^{n_U},\label{eq:sec5:updating_CBF_cons_2}
\end{align}
\end{subequations}
where $\gamma \in \mathbb{R}[x]$, $\pi \in \mathbb{R}[x]^m$, $\Omega \in \text{SOS}[x]$, and \mbox{$\Psi \in \text{SOS}[x]^{n_U}$} are unknown. As already mentioned, the SOS constraint \eqref{eq:sec5:updating_CBF_original_function} cannot be directly imposed in a standard feasibility problem in SOS programming. 

Due to the quadratic nature of $h$ and considering
the control-affine system \eqref{eq:sec2:dynamical-system}, the polynomial $h^{(k-1)}\left(x^+ \right)$ can be expressed as 
\begin{align}
    h^{(k-1)}(x^+) &= 
    \sum_{i = 1}^{m} \sum_{j = i}^{m} a^{(k-1)}_{i,j}(x)\pi_i(x)\pi_j(x) + \nonumber \\
   &\qquad ~\sum_{i = 1}^{m} b^{(k-1)}_{i}(x)\pi_i(x) + c^{(k-1)}(x),
   \label{eq:sec5:h-expressed-polynomial}
\end{align}
where $a^{(k-1)}_{i,j}, b^{(k-1)}_i, c^{(k-1)} \in \mathbb{R}[x]$, $i \in \{1,2,\hdots, m\}$, \mbox{$j \in \{i,i+1, \hdots, m\}$},
are known since $h^{(k-1)}$ is given.
To handle the bilinear terms in \eqref{eq:sec5:h-expressed-polynomial} and, consequently, in the SOS constraint \eqref{eq:sec5:updating_CBF_original_function}, we replace each $\pi_i(x)\pi_j(x)$ with an unknown $\widetilde{\pi}_{i,j} \in \mathbb{R}[x]$ and define a new polynomial $\widetilde{h}_+^{(k-1)}$. To be more precise, 
\begin{align} \label{eq:sec5:modified-DTCBF}
    \widetilde{h}^{(k-1)}_+\left(x\right) &\defeq 
    \sum_{i = 1}^{m} \sum_{j = i}^{m} a^{(k-1)}_{i,j}(x)\widetilde{\pi}_{i,j}(x) + \nonumber \\
   &\qquad \sum_{i = 1}^{m} b^{(k-1)}_{i}(x)\pi_i(x) + c^{(k-1)}(x), \hspace{-0.009cm}
\end{align}
where $\widetilde{\pi}_{i,j} \in \mathbb{R}[x]$, $i \in \{1,2,\hdots, m\}$, $j \in \{i,i+1, \hdots, m\}$.
Then, instead of \eqref{eq:sec5:updating_CBF_original_function}, we impose 
\begin{align}
    &\widetilde{h}^{(k-1)}_+(x) - h^{(k-1)}(x) + \gamma(h^{(k-1)}(x)) \nonumber \\
    & \qquad \qquad \qquad -\Lambda(x)h^{(k-1)}(x) \in \text{SOS}[x], \label{eq:sec5:updating_CBF_modified_function}
\end{align}
where $\gamma, \widetilde{\pi}_{i,j} \in \mathbb{R}[x]$ and $\Lambda \in \text{SOS}[x]$ are unknown, and in the following, we aim to introduce supplementary constraints to ensure that the satisfaction of \eqref{eq:sec5:updating_CBF_modified_function} implies that for all $x \in \mathcal{C}^{(k-1)}$, where $\mathcal{C}^{(k-1)}$ is the zero-superlevel set of $h^{(k-1)}$, it holds that \vspace{-0.2cm} 
\begin{align} \label{eq:sec5:DTCBF-aim}
    h^{(k-1)}(x^+) - h^{(k-1)}(x) + \gamma(h^{(k-1)}(x)) \geqslant 0. 
\end{align} \\[-0.6cm]
    Note that the SOS multiplier $\Omega$ in \eqref{eq:sec5:updating_CBF_original_function} may differ from the SOS multiplier $\Lambda$ in \eqref{eq:sec5:updating_CBF_modified_function}. 

Specifically, in Theorem \ref{theorem:sec5:thm1} below, we will propose supplementary constraints to ensure that 
\begin{align} \label{eq:sec5:rest-obj-const-updating-CBF-1}
    a^{(k-1)}_{i,j}(x)\bigl(\pi_i(x)\pi_j(x) - \widetilde{\pi}_{i,j}(x)\bigr) &\geqslant 0 
\end{align}
is satisfied for all $x \in \mathcal{C}^{(k-1)}$ and for all $i \in \{1,2,\hdots, m\}$ and \mbox{$j  \in \{i,i+1, \hdots, m\}$}, and as such, it holds that
\begin{align} \label{eq:sec5:obj-const-updating-CBF}
    h^{(k-1)}(x^+) \geqslant \widetilde{h}^{(k-1)}_+(x) \quad \text{for all } x \in \mathcal{C}^{(k-1)}.
\end{align}
The satisfaction of \eqref{eq:sec5:updating_CBF_modified_function} together with \eqref{eq:sec5:obj-const-updating-CBF} implies that \eqref{eq:sec5:DTCBF-aim} holds for all $x \in \mathcal{C}^{(k-1)}$. 
It is noteworthy that the condition \eqref{eq:sec5:rest-obj-const-updating-CBF-1} does not introduce conservatism since $\widetilde{\pi}_{i,j}(x)$ can be equal to $\pi_i(x)\pi_j(x)$.
To proceed, we first present Lemma \ref{lemma:sec5:det-nonnegative}, which is fundamental to our approach. 
\begin{lem}[Positive semi-definite matrices \cite{LCP-book}] \label{lemma:sec5:det-nonnegative}
 \hspace{0.05cm} A symmetric matrix $Q \defeq$ \scalebox{0.7}{$\begin{bmatrix}
         Q_{11} & Q_{12} \\
         Q_{12} & Q_{22}
     \end{bmatrix}$} $\in \mathbb{R}^{2 \times 2}$
 is positive semi-definite, denoted by $Q \succcurlyeq 0$, if and only if $Q_{11} \geqslant 0$, $Q_{22} \geqslant 0$, and $Q_{11}Q_{22} - Q_{12}^2 \geqslant 0$.
\end{lem}
To satisfy \eqref{eq:sec5:rest-obj-const-updating-CBF-1} for all \mbox{$i = j \in \{1,2, \hdots, m\}$} and for all $x \in \mathcal{C}^{(k-1)}$, we use Proposition \ref{prop:sec5:prop1} below. 
\begin{prop} \label{prop:sec5:prop1}
    Consider the DTCBF $h^{(k-1)} \in \mathbb{R}[x]$ with zero-superlevel set $\mathcal{C}^{(k-1)}$ obtained at the \mbox{$(k-1)$}-th iteration, $k \in \mathbb{N}$, for the system \eqref{eq:sec2:dynamical-system} with the control admissible set $\mathbb{U}$. Suppose there exist $\pi_i, \widetilde{\pi}_{i,i}\in \mathbb{R}[x]$ and $\sigma_{i,1}, \hdots,$ $\sigma_{i,4} \in \mathrm{SOS}[x]$, $i \in \{1,2, \hdots, m\}$, such that for all $x \in \mathbb{R}^n$, it holds that
    \begin{align} \label{eq:prop:bi-A-1}
        \begin{bmatrix}
            1 & \pi_i(x) \\
            \pi_i(x) & \mathrm{Y}_{i,1}(x)
        \end{bmatrix} \succcurlyeq 0,
    \end{align} 
    with 
    \begin{align}
        \mathrm{Y}_{i,1} \defeq \widetilde{\pi}_{i,i} - \sigma_{i,1}h^{(k-1)} + \sigma_{i,2}a^{(k-1)}_{i,i}, \nonumber
    \end{align} 
    where $a^{(k-1)}_{i,i} \in \mathbb{R}[x]$ is associated with $h^{(k-1)}(f(x)+g(x)\pi(x))$ as in \eqref{eq:sec5:h-expressed-polynomial}, and it holds that 
    \begin{align} \label{eq:prop:bi-A-2}
        &- \widetilde{\pi}_{i,i}(x) - \sigma_{i,3}(x)h^{(k-1)}(x) - \nonumber \\
        &\qquad \qquad \qquad \sigma_{i,4}(x)a^{(k-1)}_{i,i}(x) \in \mathrm{SOS}[x].
    \end{align}
    Then, for all $x \in \mathcal{C}^{(k-1)}$, it follows that 
    \begin{align} \label{eq:prop1-final}
       a^{(k-1)}_{i,i}(x)\bigl(\pi_i^2(x) - \widetilde{\pi}_{i,i}(x)\bigr) \geqslant 0.
    \end{align}
\end{prop}
\vspace{-0.7cm}
\begin{pf}
    See Appendix \ref{sec:appendix-1}.
\end{pf} 
\begin{rem}
    To impose a matrix inequality similar to \eqref{eq:prop:bi-A-1} in SOS programming, we can state that for a $Q \in \mathbb{R}[x]^{2 \times 2}$, $Q(x) \succcurlyeq 0$ holds for all $x \in \mathbb{R}^n$, if $y^{\T}Q(x)y \in \mathrm{SOS}[z]$, where $z \defeq [x^{\T} ~ y^{\T}]^{\T} \in \mathbb{R}^{n+2}$. Note that this SOS constraint preserves the linearity of decision variables.
\end{rem}
%

Moreover, we use Propositions \ref{prop:sec5:prop2} and \ref{prop:sec5:prop3} to satisfy \eqref{eq:sec5:rest-obj-const-updating-CBF-1} for all $i \in \{1,2, \hdots, m-1\}$ and $j \in \{i+1,i+2, \hdots, m\}$ and for all $x \in \mathcal{C}^{(k-1)}$.
\begin{prop}
\label{prop:sec5:prop2}
     Consider the DTCBF $h^{(k-1)} \in \mathbb{R}[x]$ with zero-superlevel set $\mathcal{C}^{(k-1)}$ obtained at the \mbox{$(k-1)$}-th iteration, $k \in \mathbb{N}$, for the system \eqref{eq:sec2:dynamical-system} with the control admissible set $\mathbb{U}$. Suppose there exist $\pi_1, \hdots, \pi_m, \widetilde{\pi}_{i,j},$ $\Theta_{i,j,1},\Theta_{i,j,2},\Theta_{i,j,3} \in \mathbb{R}[x]$ and $\xi_{i,j,1}, \hdots, \xi_{i,j,8} \in \mathrm{SOS}[x]$, $i \in \{1,2, \hdots, m - 1\}$, $j \in \{i+1,i+2, \hdots, m\}$, such that for all $x \in \mathbb{R}^n$, it holds that 
    \begin{align}
        \left[\begin{matrix}
         1 & \pi_i(x) \\
         \pi_i(x) &  \Upsilon_{i,j,1}(x)
        \end{matrix}
        \right] &\succcurlyeq 0, \label{eq:prop:bi-B-1} \\
         \left[\begin{matrix}
         1 & \pi_j(x) \\
         \pi_j(x) &  \Upsilon_{i,j,2}(x)
        \end{matrix}
        \right] &\succcurlyeq 0, \label{eq:prop:bi-B-2}
        \end{align}
        with
        \begin{align}
            \Upsilon_{i,j,1} &\defeq \Theta_{i,j,1} - \xi_{i,j,1}h^{(k-1)} + \xi_{i,j,2}a^{(k-1)}_{i,j}, \nonumber \\
            \Upsilon_{i,j,2} &\defeq \Theta_{i,j,2} - \xi_{i,j,3}h^{(k-1)} + \xi_{i,j,4}a^{(k-1)}_{i,j}, \nonumber
        \end{align}
    where $a^{(k-1)}_{i,j} \in \mathbb{R}[x]$ is associated with $h^{(k-1)}(f(x)+g(x)\pi(x))$ as in \eqref{eq:sec5:h-expressed-polynomial}, and it holds that
    \begin{align}
        &\Theta_{i,j,3}(x) - \xi_{i,j,5}(x)h^{(k-1)}(x) + \nonumber \\
        &\hspace{3.132cm} \xi_{i,j,6}(x)a^{(k-1)}_{i,j}(x) \in \mathrm{SOS}[x], \label{eq:prop:bi-B-3} \\
        &2\widetilde{\pi}_{i,j}(x) - (\Theta_{i,j,1}(x) + \Theta_{i,j,2}(x) + \Theta_{i,j,3}(x)) - \nonumber \\
        &~\xi_{i,j,7}(x)h^{(k-1)}(x) + \xi_{i,j,8}(x)a^{(k-1)}_{i,j}(x) \in \mathrm{SOS}[x]. \label{eq:prop:bi-B-4}
    \end{align} 
    Then, for all $x \in \mathcal{C}^{(k-1)}$ such that $a^{(k-1)}_{i,j}(x) \leqslant 0$, it follows that 
    \begin{align}
        \widetilde{\pi}_{i,j}(x) - \pi_i(x)\pi_j(x) \geqslant 0. \nonumber
    \end{align}
\end{prop}
\begin{pf}
    See Apendix \ref{sec:appendix-2}.
\end{pf}
\begin{prop} \label{prop:sec5:prop3} 
    Consider the DTCBF \mbox{$h^{(k-1)} \in \mathbb{R}[x]$} with zero-superlevel set $\mathcal{C}^{(k-1)}$ obtained at the \mbox{$(k-1)$}-th iteration, $k \in \mathbb{N}$, for the system \eqref{eq:sec2:dynamical-system} with the control admissible set $\mathbb{U}$. Suppose there exist $\pi_1, \hdots, \pi_m, \widetilde{\pi}_{i,j},$ $\Delta_{i,j,1}, \Delta_{i,j,2}, \Delta_{i,j,3} \in \mathbb{R}[x]$ and $\eta_{i,j,1}, \hdots,$ $\eta_{i,j,8} \in \mathrm{SOS}[x]$, $i \in \{1,2,\hdots, m-1\}$, $j \in\{i+1,i+2, \hdots, m\}$, such that for all $x \in \mathbb{R}^n$, it holds that 
    \begin{align} 
        \left[\begin{matrix}
         1 & \pi_i(x) \\
         \pi_i(x) &  \mathrm{T}_{i,j,1}(x)
        \end{matrix}
        \right] &\succcurlyeq 0, \label{eq:prop:bi-B-5}  \\
         \left[\begin{matrix}
         1 & \pi_j(x) \\
         \pi_j(x) &  \mathrm{T}_{i,j,2}(x)
        \end{matrix}
        \right] &\succcurlyeq 0, \label{eq:prop:bi-B-6}
    \end{align}
        with 
    \begin{align}
        \mathrm{T}_{i,j,1} &\defeq \Delta_{i,j,1} - \eta_{i,j,1}h^{(k-1)} -\eta_{i,j,2}a^{(k-1)}_{i,j}, \nonumber\\
        \mathrm{T}_{i,j,2} &\defeq \Delta_{i,j,2} - \eta_{i,j,3}h^{(k-1)} -\eta_{i,j,4}a^{(k-1)}_{i,j}, \nonumber
    \end{align}
    where $a^{(k-1)}_{i,j} \in \mathbb{R}[x]$ is associated with $h^{(k-1)}(f(x)+g(x)\pi(x))$ as in \eqref{eq:sec5:h-expressed-polynomial}, and it holds that
    \begin{align}
         &- \Delta_{i,j,3}(x) - \eta_{i,j,5}(x)h^{(k-1)}(x) - \nonumber \\
        &\hspace{3.55cm} \eta_{i,j,6}(x)a^{(k-1)}_{i,j} \in \mathrm{SOS}[x], \label{eq:prop:bi-B-7} \\
        &(\Delta_{i,j,3}(x) - \Delta_{i,j,1}(x) - \Delta_{i,j,2}(x)) - 2\widetilde{\pi}_{i,j}(x) - \nonumber \\
        &\eta_{i,j,7}(x)h^{(k-1)}(x) - \eta_{i,j,8}(x)a^{(k-1)}_{i,j}(x)\in \mathrm{SOS}[x]. \label{eq:prop:bi-B-8}
    \end{align}
    Then, for all \mbox{$x \in \mathcal{C}^{(k-1)}$} such that $a_{i,j}(x) \geqslant 0$, it follows that
    \begin{align}
     \pi_i(x)\pi_j(x) - \widetilde{\pi}_{i,j}(x) \geqslant 0. \nonumber
    \end{align}
\end{prop}
\begin{pf}
    Similar to the proof of Proposition \ref{prop:sec5:prop2}. 
\end{pf}
Combining all of the above auxiliary results, we can formulate SOS constraints where decision variables appear linearly, ensuring that \eqref{eq:sec5:rest-obj-const-updating-CBF-1} is satisfied. This allows us to update the control policy $\pi$ and $\gamma$ at the $k$-th iteration, given the DTCBF $h^{(k-1)}$ obtained at the $(k-1)$-th iteration.
\begin{thm} \label{theorem:sec5:thm1}
    Consider the DTCBF $h^{(k-1)} \in \mathbb{R}[x]$ with zero-superlevel set $\mathcal{C}^{(k-1)}$ obtained at the $(k − 1)$-th iteration, $k \in \mathbb{N}$, for the system \eqref{eq:sec2:dynamical-system} with the control admissible set $\mathbb{U}$. Assume there exist $\pi_p, \widetilde{\pi}_{p,q} \in \mathbb{R}[x]$, $p \in \{1,2,\hdots, m\}$, $q \in \{p,p+1, \hdots, m\}$, $\Lambda \in \mathrm{SOS}[x]$, and $\Psi \in \mathrm{SOS}[x]^{n_U}$ such that \eqref{eq:sec5:updating_CBF_cons_2} is satisfied and the modified DTCBF $\widetilde{h}^{(k-1)}_+$ constructed as in \eqref{eq:sec5:modified-DTCBF} satisfies \eqref{eq:sec5:updating_CBF_modified_function}. Moreover, assume that there exist $\sigma_{p,1}, \hdots, \sigma_{p,4} \in \mathbb{R}[x]$ such that \eqref{eq:prop:bi-A-1} and \eqref{eq:prop:bi-A-2} are satisfied for all $p\in \{1,2, \hdots, m\}$, and there exist $\xi_{i,j,1}, \hdots, \xi_{i,j, 8}, \eta_{i,j,1}, \hdots, \eta_{i,j,8} \in \mathrm{SOS}[x]$ as well as $\Theta_{i,j,1},\Theta_{i,j,2}, \Theta_{i,j,3},\Delta_{i,j,1}, \Delta_{i,j,2}, \Delta_{i,j,3} \in \mathbb{R}[x]$ such that \eqref{eq:prop:bi-B-1}\textendash\eqref{eq:prop:bi-B-8} are satisfied for all $i \in \{1,2, \hdots, m-1\}$ and $j\in \{i+1,i+2, \hdots, m\}$. Then, it follows that for all $x \in \mathcal{C}^{(k-1)}$, it holds that
    \begin{align} 
        &h^{(k-1)}(f(x) + g(x)\pi(x)) - \nonumber \\
        &\qquad \qquad h^{(k-1)}(x) + \gamma(h^{(k-1)}(x)) \geqslant 0, \label{eq:thm1:eq1}
    \end{align}
    where $\pi(x) \defeq [\pi_1(x) ~ \hdots ~ \pi_m(x)]^{\T}$, and $\pi$ is admissible for all $x \in \mathcal{C}^{(k-1)}$ with respect to $\mathbb{U}$.
\end{thm}
\vspace{-0.3cm}
\begin{pf}
    Imposing \eqref{eq:sec5:updating_CBF_cons_2} implies that $\pi$ is admissible for all $x \in \mathcal{C}^{(k-1)}$ with respect to $\mathbb{U}$ based on Lemma \ref{lemma:sec4:S-procedure}. \\
    By Proposition \ref{prop:sec5:prop1}, the satisfaction of \eqref{eq:prop:bi-A-1} and \eqref{eq:prop:bi-A-2} implies that for all $i \in \{1, \hdots, m\}$ and \mbox{for all $x \in \mathcal{C}^{(k-1)}$}, \eqref{eq:prop1-final} holds.
    Besides that, based on Propositions \mbox{\ref{prop:sec5:prop2} and \ref{prop:sec5:prop3}}, the satisfaction of \eqref{eq:prop:bi-B-1}\textendash \eqref{eq:prop:bi-B-8} implies that for all \mbox{$i \in \{1,2, \hdots, m-1\}$}, $j \in \{i+1,i+2, \hdots, m\}$, and for all \mbox{$x \in \mathcal{C}^{(k-1)}$}, 
    \begin{align} \label{eq:proof:thm1:eq2}
        a^{(k-1)}_{i,j}(x)(\pi_i(x)\pi_j(x) - \widetilde{\pi}_{i,j}(x)) \geqslant 0.
    \end{align}
    By virtue of \eqref{eq:prop1-final} and \eqref{eq:proof:thm1:eq2}, we can conclude that for all \mbox{$x \in \mathcal{C}^{(k-1)}$}, 
    \begin{align} \label{eq:proof:thm1:eq3}
        h^{(k-1)}(f(x)+g(x)\pi(x)) \geqslant \widetilde{h}^{(k-1)}_+(x).
    \end{align}
     As a result, satisfying \eqref{eq:proof:thm1:eq3} together with \eqref{eq:sec5:updating_CBF_modified_function} implies that \eqref{eq:thm1:eq1} holds.\qed
\end{pf}


\subsection{Algorithm Overview}\label{sec:sec5.2}
The iterative algorithm to synthesize a DTCBF-triple $(h, \gamma, \pi)$ for the system \eqref{eq:sec2:dynamical-system} with the control admissibe set $\mathbb{U}$ is proposed as follows: \vspace{0.2cm}
\begin{enumerate}[listparindent=\parindent, label=\textbf{Step \arabic*)}, ref=\arabic*, wide=0pt] \addtocounter{enumi}{-1}
    \item \textbf{Initialization:}  \label{step:step0} We parameterize $h$, $\gamma$, and $\pi$ as in \eqref{eq:sec4:parameterized-DTCBF}, \eqref{eq:sec4:parameterized-gamma}, and \eqref{eq:sec4:parameterized-control-policy}, respectively. We also set the degree of auxiliary polynomials  $\widetilde{\pi}_{p,q} \in \mathbb{R}[x]$, $p \in \{1,2, \hdots, m\}$, \mbox{$q \in \{p,p+1, \hdots, m\}$}. We set the degree of SOS multipliers $\Lambda, \Omega, \Phi, \sigma_{p,1}, \hdots, \sigma_{p,4} \in \mathrm{SOS}[x]$, \mbox{$p \in \{1,2, \hdots, m\}$}, $\Psi \hspace{-0.02cm} \in \hspace{-0.02cm} \mathrm{SOS}[x]^{n_U}$, $\xi_{i,j,1}, \hdots, \xi_{i,j,8}, \eta_{i,j,1},\hdots, \eta_{i,j, 8} \in \mathrm{SOS}[x]$, and $\Theta_{i,j,1}, \Theta_{i,j,2}, \Theta_{i,j,3}, \Delta_{i,j,1},\Delta_{i,j,2},\Delta_{i,j,3} \in \mathbb{R}[x]$, \mbox{$i \in \{1,2, \hdots, m-1\}$} and $j \in \{i+1,i+2, \hdots, m\}$. 
    
We initialize the DTCBF $h^{(0)}$, which must be a valid quadratic DTCBF, though it may be quite conservative. For example, if the origin is a stable fixed point, a small region around it can be considered.
We set the iteration counter to $k\defeq1$.
\vspace{0.15cm}
\item \textbf{Updating the Control Policy:} \label{step:step1} Based on Theorem \ref{theorem:sec5:thm1} and given $h^{(k-1)}$ obtained at the $(k-1)$-th iteration, we construct a feasibility problem in SOS programming to update $\gamma$ and its friend $\pi$ as \vspace{-0.2cm}
\begin{subequations} \label{eq:step1:opti1}
\begin{align}
    \mathrm{find} ~&~ \gamma, \pi_p, \widetilde{\pi}_{p,q}, \Lambda, \Psi, \sigma_{p,1}, \hdots, \sigma_{p,4}, \nonumber \\
    &~\xi_{i,j,1}, \hdots, \xi_{i,j,8}, \eta_{i,j,1}, \hdots, \eta_{i,j,8}, \nonumber \\
    &~\Theta_{i,j,1}, \Theta_{i,j,2}, \Theta_{i,j,3}, \Delta_{i,j,1}, \Delta_{i,j,2}, \Delta_{i,j,3} \\
    \mathrm{s.t.} ~&~ \eqref{eq:sec5:updating_CBF_cons_2}, \eqref{eq:sec5:updating_CBF_modified_function}, \eqref{eq:prop:bi-A-1}, \eqref{eq:prop:bi-A-2}, \eqref{eq:prop:bi-B-1}\textrm{--}\eqref{eq:prop:bi-B-8}, \nonumber \\
    & \gamma, \pi_p, \widetilde{\pi}_{p,q},\Theta_{i,j,1}, \Theta_{i,j,2}, \Theta_{i,j,3}, \nonumber\\
    & \hspace{1.55cm} \Delta_{i,j,1}, \Delta_{i,j,2}, \Delta_{i,j,3} \in \mathbb{R}[x], \\
    &~ \Lambda, \sigma_{p,1}, \hdots, \sigma_{p,4}, \xi_{i,j,1}, \hdots, \xi_{i,j,8}, \nonumber \\
    & \hspace{2.172cm} \eta_{i,j,1}, \hdots, \eta_{i,j,8} \in \mathrm{SOS}[x], \\
    &\hspace{4.1cm} \Psi \in \mathrm{SOS}[x]^{n_U},
\end{align}
\end{subequations}
where $p\in \{1,2, \hdots, m\}$, $q \in \{p,p+1, \hdots, m\}$, \mbox{$i \in \{1,2, \hdots, m-1\}$}, and \mbox{$j \in \{i+1,i+2, \hdots, m\}$}. \vspace{0.3cm}

Using $\gamma^{(k)}$ and $\pi^{(k)}$ obtained in \eqref{eq:step1:opti1} at the current iteration $k$, we aim to find \mbox{a SOS multiplier $\Omega \in \mathrm{SOS}[x]$} such that \eqref{eq:sec5:updating_CBF_original_function} holds. This multiplier is
essential for updating the DTCBF in the subsequent step. To find it, we construct a feasibility problem in SOS programming as \vspace{-0.2cm}
\begin{subequations} \label{eq:step1:opti2}
\begin{align}
    \mathrm{find} ~&~ \Omega \\
    \mathrm{s.t.} ~&~ h^{(k-1)}(f(x) + g(x)\pi^{(k)}(x)) - h^{(k-1)}(x) + \nonumber \\
    & \gamma^{(k)}(h^{(k-1)}(x)) - \Omega(x) h^{(k-1)}(x) \in \mathrm{SOS}[x], \\
    &\hspace{4.64cm} \Omega \in \mathrm{SOS}[x].
\end{align}
\end{subequations}
\item \textbf{Updating the DTCBF:} \label{step:step2}
In this step, we obtain the updated DTCBF $h^{(k)}$ by solving the feasibility problem \vspace{-0.2cm} 
\begin{subequations} \label{eq:step2:opti1}
\begin{align}
    &\hspace{-5cm} \mathrm{find} ~~ h, \Phi, \Xi \\
    \mathrm{s.t.} ~~ \eqref{eq:sec5:updating_CBF_modified_function}, \eqref{eq:prop:bi-A-1}, \eqref{eq:prop:bi-A-2}, \eqref{eq:prop:bi-B-1}\textrm{--}\eqref{eq:prop:bi-B-8}, \hspace{0.5cm}  \nonumber \\
    h(f(x) + g(x) \pi^{(k)}(x)) - h(x&)  +  \nonumber \\
    \qquad \gamma^{(k)}(h(x)) - \Omega^{(k)}(x)h(x) &\in \mathrm{SOS}[x], \label{eq:step2:opti1-const-1} \\
    M\pi^{(k)}(x) + d - \Psi^{(k)}(x)h(x) &\in \mathrm{SOS}[x]^{n_U}, \label{eq:step2:opti1-const-3} \\
    -h(x) - \epsilon + \Phi(x)s(x) &\in \mathrm{SOS}[x], \label{eq:step2:opti1-const-2} \\
    h(x) - \delta - \Xi(x) h^{(k-1)}(x) &\in \mathrm{SOS}[x], \label{eq:step2:opti1-const-4} \\
    \Phi, \Xi &\in \mathrm{SOS}[x], \label{eq:step2:opti1-const-5} \\
    h
    &\in \mathbb{R}[x].
\end{align}
\end{subequations}
Here, \mbox{$\epsilon, \delta \in \mathbb{R}_{>0}$} are predefined and sufficiently small, $M \in \mathbb{R}^{n_U \times m}$ and $d \in \mathbb{R}^{n_U}$ are associated with the control admissible set $\mathbb{U}$ in \eqref{eq4:sec4:control-admissible-assumption}, and the mapping $s$ is associated with the safe set $\mathcal{S}$ in \eqref{eq:sec2:safe-set}. Additionally, $\pi^{(k)}$, $\gamma^{(k)}$, and $\Psi^{(k)}$ are the solutions of the feasibility \mbox{problem \eqref{eq:step1:opti1}}, and $\Omega^{(k)}$ is the solution of the \mbox{problem \eqref{eq:step1:opti2}} at the current iteration $k$. Moreover, $h^{(k-1)}$ is the solution of \eqref{eq:step2:opti1} at the previous iteration.
\vspace{0.2cm}
\begin{rem} \label{rem:sec5:lossless}
    We impose the SOS constraints \eqref{eq:sec5:updating_CBF_modified_function}, \eqref{eq:prop:bi-A-1}, \eqref{eq:prop:bi-A-2}, \eqref{eq:prop:bi-B-1}\textrm{--}\eqref{eq:prop:bi-B-8} in the feasibility problem \eqref{eq:step2:opti1} to ensure that the problem \eqref{eq:step1:opti1} is feasible in the next iteration.
    In these constraints, the SOS multipliers are fixed and obtained in \eqref{eq:step1:opti1} at the current iteration $k$, 
    while \mbox{$h, a_{i,j} \in \mathbb{R}[x]$}, $i \in \{1,2, \hdots, m\}$, \mbox{$j \in \{i,i+1, \hdots, m\}$}, are unknown. However, we do not restate the constraints to avoid redundancy.
\end{rem}
\vspace{0.2cm}
\begin{rem}
    One may consider maximizing $\delta$ as an objective in the problem \eqref{eq:step2:opti1} to maximize the increase in the size of the zero-superlevel set relative to that obtained in the previous iteration, $\mathcal{C}^{(k-1)}$.
\end{rem}
\vspace{0.2cm}
If the problem \eqref{eq:step2:opti1} is feasible, we proceed to \mbox{Step \ref{step:step1}} and set \mbox{$k \defeq k +1$}. If it is infeasible due to violation of \eqref{eq:step2:opti1-const-4}, indicating that the size of the zero-superlevel set of the DTCBF can no longer be enlarged, the algorithm is terminated, and we report $(h^{(k-1)}, \gamma^{(k-1)},\pi^{(k-1)})$ as the synthesized DTCBF-triple.

\end{enumerate}

\begin{thm}
    The proposed algorithm is lossless, meaning that once the conditions in Step \ref{step:step1} are satisfied, i.e., the control policy $\pi^{(1)}$ and $\gamma^{(1)}$ are found for the initial DTCBF $h^{(0)}$ corresponding to the system \eqref{eq:sec2:dynamical-system} with the control admissible set $\mathbb{U}$, the satisfaction of all conditions is guaranteed, except \eqref{eq:step2:opti1-const-4} which indicates that the size of the zero-superlevel set of $h$ can no longer be enlarged. Moreover, it is guaranteed that $(h^{(l)}, \gamma^{(l)}, \pi^{(l)})$, for all $l \in \{1,2,\hdots, k-1\}$, forms a DTCBF-triple for the \mbox{system \eqref{eq:sec2:dynamical-system}} with $\mathbb{U}$, and that \mbox{$\mathcal{C}^{(l-1)} \subset \mathcal{C}^{(l)} \subseteq \mathcal{S}$} holds. 
\end{thm}
\vspace{-0.5cm}
\begin{pf}
    Imposing the conditions \eqref{eq:sec5:updating_CBF_modified_function}, \eqref{eq:prop:bi-A-1}, \eqref{eq:prop:bi-A-2}, \eqref{eq:prop:bi-B-1}\textrm{--}\eqref{eq:prop:bi-B-8}, and \eqref{eq:step2:opti1-const-3} in the feasibility problem \eqref{eq:step2:opti1} in Step \ref{step:step2} at the $k$-th iteration, $k \in \mathbb{N}$, ensures that
    all conditions of the problem \eqref{eq:step1:opti1} in Step \ref{step:step1} at the $(k+1)$-th iteration can be satisfied, even in the worst case, with $\pi^{(k)}$ and $\gamma^{(k)}$ obtained at the $k$-th iteration. Subsequently, the feasibility of the \mbox{problem \eqref{eq:step1:opti1}} at the $k$-th iteration, $k \in \mathbb{N}$, guarantees that all conditions of \mbox{Step \ref{step:step2}}, except \eqref{eq:step2:opti1-const-4}, at this iteration can be satisfied with $h^{(k-1)}$. Therefore, the only condition that may halt the algorithm is if the size of the zero-superlevel set of $h$ can no longer be enlarged.
    Additionally, based on Lemma \ref{lemma:sec4:S-procedure}, the satisfaction of \eqref{eq:step2:opti1-const-1} and \eqref{eq:step2:opti1-const-3} implies that $(h^{(l)}, \gamma^{(l)}, \pi^{(l)})$, for all $l \in \{1,2, \hdots, k-1\}$, is a DTCBF-triple. Lastly, we impose \eqref{eq:step2:opti1-const-2} and \eqref{eq:step2:opti1-const-4} to ensure that $\mathcal{C}^{(l-1)} \subset \mathcal{C}^{(l)} \subseteq \mathcal{S}$. \qed
\end{pf}

\begin{rem}
    Note that the feasibility of \eqref{eq:step1:opti1} at the first iteration ($k=1$) is not generally guaranteed due to the conservatism introduced both by the supplementary conditions used to handle the bilinear terms and by the use of the generalized S-procedure lemma (Lemma \ref{lemma:sec4:S-procedure}). If \eqref{eq:step1:opti1} is initially infeasible, one may consider selecting a different initial DTCBF $h^{(0)}$. It is worth mentioning that an initial DTCBF whose zero-superlevel set lies entirely within a circular region may represent a best-case scenario. In this case, the only conservatism arises from the use of the generalized S-procedure lemma, as we only need the condition \eqref{eq:prop:bi-A-1} to handle the bilinear terms (see the proof of Proposition \ref{prop:sec5:prop1} for more details). 
\end{rem}

\section{Proposed Synthesis Method for DTCBFs with Higher Polynomial Degrees} \label{sec:sec6}
This section explains two approaches by which the method proposed in Section \ref{sec:sec5} extends to the synthesis of polynomial DTCBFs of degree higher than two for the system \eqref{eq:sec2:dynamical-system} with the control admissible set~$\mathbb{U}$. 
\vspace{-0.18cm}
\subsection{Multiple Uses of Propositions \ref{prop:sec5:prop1}--\ref{prop:sec5:prop3}} \label{sec:sec8-1}
\vspace{-0.18cm}
The first approach is to apply the conditions in Propositions \ref{prop:sec5:prop1}--\ref{prop:sec5:prop3} multiple times to handle all the bilinear terms. To do so, we first redefine the control input $u$ in the system \eqref{eq:sec2:dynamical-system} by introducing a shifted control input $\widetilde{u} \defeq u + c$. Here, $c \in \mathbb{R}^m$ is chosen to be sufficiently large so that $\widetilde{u}_i \geqslant 0$ for all $i \in \{1, \hdots, m\}$ and for all $u \in \mathbb{U}$, assuming $\mathbb{U}$ is bounded. This transforms the system \eqref{eq:sec2:dynamical-system} into \vspace{-0.1cm}
\begin{align}
    x^+ = f(x) + g(x)u = \widetilde{f}(x) + g(x)\widetilde{u},
\end{align} \\ [-0.4cm]
where $f, \widetilde{f} \in \mathbb{R}[x]^n$, $g \in \mathbb{R}[x]^{n \times m}$, and $\widetilde{u} \in \widetilde{\mathbb{U}} \subset \mathbb{R}^{m}_{\geqslant 0}$.
Rather than synthesizing a control policy $\pi$ for $u$, we aim to synthesize a control policy $\mu$ for $\widetilde{u}$, which can subsequently be shifted back to $\pi(x) = \mu(x) - c$.

In $h(\widetilde{f}(x) + g(x)\mu(x))$, the terms $\mu_i$, $i \in \{1,2, \hdots, m\}$, where $\mu(x) \defeq [\mu_1(x) ~\hdots ~\mu_m(x)]^{\T}$, may appear in products of degrees greater than two, such as $\mu_1^3(x)$, $\mu_1^4(x)$, $\mu_1^2(x)\mu_{2}(x)$, and so on. The products of degree at most two can be addressed as discussed in \mbox{Section \ref{sec:sec5}}.
This section explains how the proposed method extends to address these higher-degree products in \mbox{Step \ref{step:step1}}, while the other steps remain unchanged. We explain it for a specific bilinear term, namely $\mu_i^3$, \mbox{$i \in \{1,2, \hdots, m\}$}. 
All other bilinear terms can be handled in a similar manner. Similar to \mbox{Section \ref{sec:sec5.1}}, to handle $\mu_i^3$, we introduce an unknown polynomial \mbox{$\widetilde{\mu}_{i,i,i} \in \mathbb{R}[x]$} and aim to impose SOS constraints where the decision variables appear linearly so that \vspace{-0.05cm}
\begin{align} \label{eq:sec8:aim}
    a^{(k-1)}_{i,i,i}(x)(\mu^3_i(x) - \widetilde{\mu}_{i,i,i}(x)) \geqslant 0, ~ i \in \{1, \hdots, m\},
\end{align} \\ [-0.4cm]
where $a_{i,i,i}^{(k-1)} \in \mathbb{R}[x]$ is known and associated with $h^{(k-1)}(\widetilde{f}(x)+g(x)\mu(x))$, as written in a form similar to \eqref{eq:sec5:h-expressed-polynomial}. To satisfy \eqref{eq:sec8:aim}, we first define a new polynomial \mbox{$\widetilde{\mu}_{i,i} \in \mathbb{R}[x]$} and impose SOS constraints based on Proposition \ref{prop:sec5:prop1}, such that for all $x \in \mathcal{C}^{(k-1)}$ satisfying $a_{i,i,i}^{(k-1)}(x) \leqslant 0$, it holds that 
\begin{align} \label{eq:sec8:eq1}
    \widetilde{\mu}_{i,i}(x) \geqslant \mu_i^2(x).
\end{align}
Then, using Proposition \ref{prop:sec5:prop2}, we impose SOS constraints, ensuring that for all $x \in \mathcal{C}^{(k-1)}$ such that $a_{i,i,i}^{(k-1)}(x) \leqslant 0$, it holds that
\begin{align} \label{eq:sec8:eq2}
    \widetilde{\mu}_{i,i,i}(x) \geqslant \widetilde{\mu}_{i,i}(x) \mu_i(x).
\end{align}
The satisfaction of \eqref{eq:sec8:eq1} and \eqref{eq:sec8:eq2}, together with the condition that $\mu_i(x) \geqslant 0$ holds for all $i \in \{1,2, \hdots, m\}$ and for all $x \in \mathcal{C}^{(k-1)}$, implies that for all $x \in \mathcal{C}^{(k-1)}$ such that $a_{i,i,i}^{(k-1)}(x) \leqslant 0$, it holds that 
\begin{align} \label{eq:sec8:eq3}
    \widetilde{\mu}_{i,i,i}(x) \geqslant \mu_i^3(x).
\end{align}
Then, similar to the condition \eqref{eq:prop:bi-A-2} in Proposition \ref{prop:sec5:prop1}, by imposing 
\begin{align}
        &- \widetilde{\mu}_{i,i,i}(x) - \widetilde{\sigma}_{i,1}(x)h^{(k-1)}(x) - \nonumber \\
        &\qquad \qquad \qquad \widetilde{\sigma}_{i,2}(x)a^{(k-1)}_{i,i,i}(x) \in \mathrm{SOS}[x],
\end{align}
where $\widetilde{\sigma}_{i,1}, \widetilde{\sigma}_{i,2} \in \text{SOS}[x]$, we can ensure that for all \mbox{$x \in \mathcal{C}^{(k-1)}$} such that \mbox{$a^{(k-1)}_{i,i,i} \geqslant 0$}, it holds that 
\begin{align} \label{eq:sec8:eq4} 
\mu_i^3(x) \geqslant 0 \geqslant \widetilde{\mu}_{i,i,i}(x).
\end{align}
In conclusion, the satisfaction of \eqref{eq:sec8:eq3} and \eqref{eq:sec8:eq4} implies that \eqref{eq:sec8:aim} holds.

\subsection{Fixed Control Policy}
The second approach involves synthesizing a quadratic DTCBF and a corresponding control policy for the system \eqref{eq:sec2:dynamical-system} with the control admissible set $\mathbb{U}$, using the method proposed in Section \ref{sec:sec5}, and progressively enlarging its zero-superlevel set until no further enlargement is feasible, i.e., until the problem \eqref{eq:step2:opti1} becomes infeasible. Subsequently, a higher-degree DTCBF can be synthesized using only Step \ref{step:step2} of the algorithm, given the control policy synthesized for the quadratic DTCBF, as the bilinearity issue does not arise when the control policy is fixed. This approach is similar to backup policy methods \cite{Gurriet2020, chen2021} for synthesizing continuous-time CBFs and is computationally more efficient than the first one (Section \ref{sec:sec8-1}). However, when the control admissible set $\mathbb{U}$ is tight, it may not significantly enlarge the zero-superlevel set of the resulting DTCBF.

\section{Numerical Case Study} \label{sec:sec7}
\vspace{-0.2cm}
To show the effectiveness of the proposed synthesis method, we apply it to a linear and a nonlinear system. For the linear system, we parameterize the DTCBF as a degree-four polynomial (based on Section \ref{sec:sec6}), and for the nonlinear system, we parameterize the DTCBF as a quadratic polynomial (based on Section \ref{sec:sec5}). The numerical results are obtained on a laptop with an Intel i7-12700H processor and Matlab 2024a, and the SOS programming problems are solved
using SOSTOOLS \cite{sostools}. 
\subsection{Linear System} 
\begin{figure}[!ht]\centering 
\includegraphics[width=0.42\textwidth]{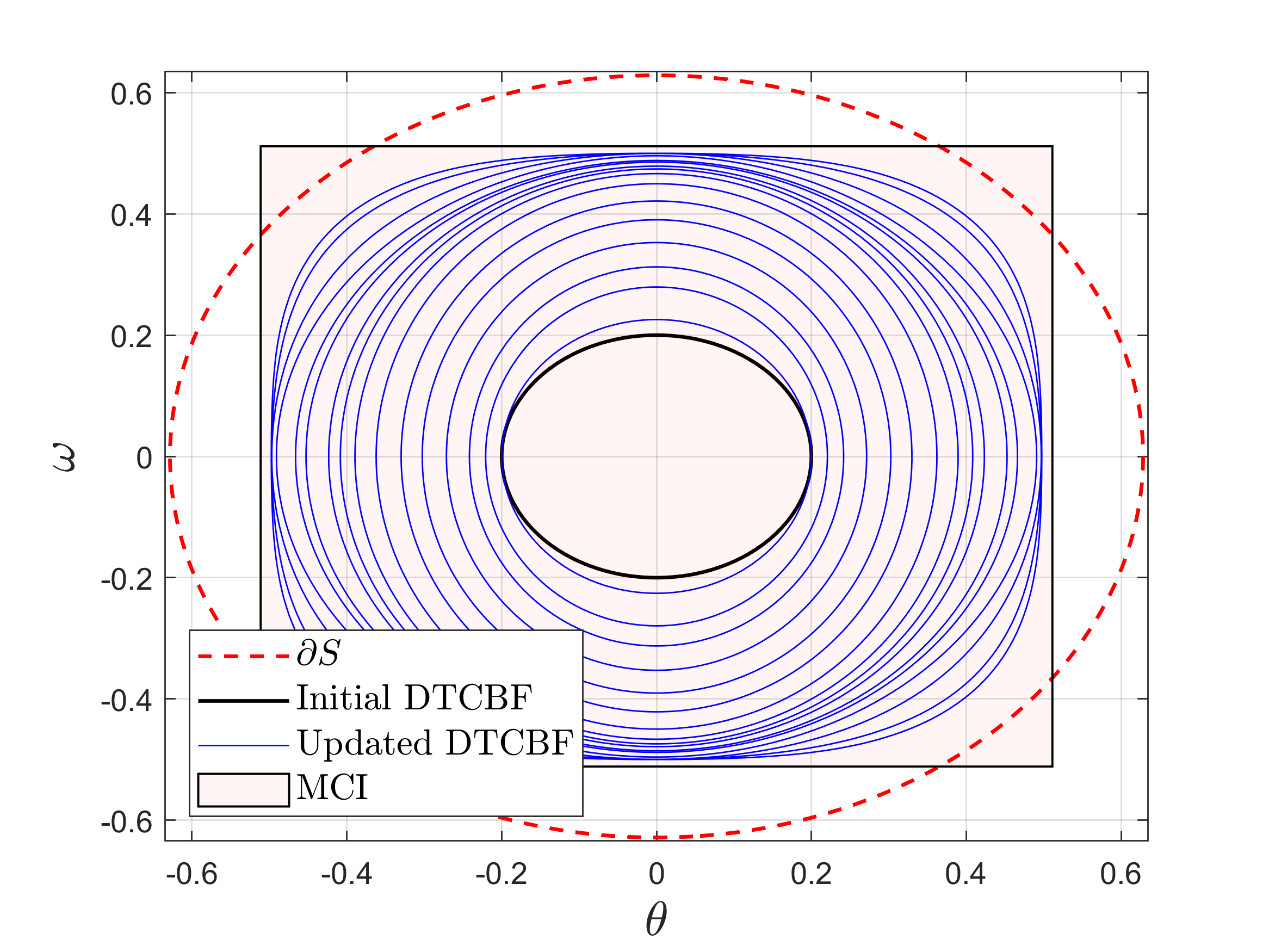}
    \caption{The proposed synthesis method is applied to the system \eqref{eq:sec6:linear-system} with the control admissible set $\mathbb{U}$ and the \mbox{safe set \eqref{eq:sec6:cart-pole-safe-set}} (the interior and boundary of the red circle) to iteratively synthesize a DTCBF-triple, starting from the initial DTCBF \eqref{eq:sec6:cart-pole-initial-guess} (black circle). Each of these blue curves represents the updated DTCBF after each iteration. The pale red square represents the maximal polyhedral controlled invariant set obtained using the MPT3 toolbox \cite{MPT3}.}
    \label{fig:sec6:example-cart-pole}
\end{figure}

Consider the linearized and discretized cart-pole system discussed in \cite{wu2023neural} as
\begin{align} \label{eq:sec6:linear-system}
    \begin{bmatrix}
        x_c^+ \\ v_c^+ \\ \theta^+ \\ \omega^+
    \end{bmatrix} = \begin{bmatrix}
        0 & 1 & 0 & 0 \\
        0 & 0 & -0.98 & 0 \\
        0 & 0 & 0 & 1 \\
        0 & 0 & 10.78 & 0
    \end{bmatrix}\begin{bmatrix}
        x_c \\ v_c \\ \theta \\ \omega
    \end{bmatrix} + \begin{bmatrix}
        0 \\ 1 \\ 0 \\ -1
    \end{bmatrix}u,
\end{align}
where $x_c, v_c \in \mathbb{R}$ are the horizontal position and velocity of the cart, respectively, $\theta \in \mathbb{R}$ is the angle of the pole from the upward vertical direction, and $\omega \in \mathbb{R}$ is its angular velocity. Additionally, $x \defeq [x_c ~ v_c ~ \theta ~ \omega]^{\T}$, and \mbox{$u \in \mathbb{U} \defeq [-5, ~ 5]$} is the horizontal force applied to the cart, which is considered as the control input to the system \eqref{eq:sec6:linear-system}. The safe set $\mathcal{S}$ is defined as the disk 
\begin{align} \label{eq:sec6:cart-pole-safe-set}
    \mathcal{S} \defeq \{ x \in \mathbb{R}^4 \mid \theta^2 + \omega^2  \leqslant (\pi/5)^2 \}.
\end{align} 
We aim to apply our proposed method to synthesize a DTCBF-triple $(h, \gamma, \pi)$ according to the objectives of \mbox{Problem \ref{prob:sec3:synthesis}} for the system \eqref{eq:sec6:linear-system} with the control admissible set $\mathbb{U}$ and the safe set \eqref{eq:sec6:cart-pole-safe-set}. For the initialization step of the proposed synthesis method (Step \ref{step:step0}), a relatively small circle, defined as 
\begin{align} \label{eq:sec6:cart-pole-initial-guess}
    h^{(0)}(x) \defeq - \theta^2 - \omega^2 + 0.04,
\end{align} 
is considered as an initial DTCBF. Moreover, we use a vector of monomials up to degrees four and three for the DTCBF $h$ and the control policy $\pi$, parameterized as in \eqref{eq:sec4:parameterized-DTCBF-general} and \eqref{eq:sec4:parameterized-control-policy}, respectively. Regarding $\gamma \in \mathcal{K}^{\leqslant \mathrm{id}}_{\infty}$, we parameterize it as in \eqref{eq:sec4:parameterized-gamma}, and assume that we aim to find $\gamma_0$ as close as possible to $0.8$, which is introduced as an objective in the SOS programming problem \eqref{eq:step1:opti1}.
After around 15 iterations, each taking an average time of $\unit[20]{s}$ (total time about $\unit[5]{min}$), the DTCBF-triple $(h, \gamma, \pi)$ is synthesized as 
\begin{align}
    h(x) &\defeq -3.910\omega^4 - 4.261\omega^2\theta^2 - 4.101\theta^4 + 0.860\omega^2 + \nonumber \\
    &\qquad \quad  0.918\theta^2 + 0.027, \nonumber \\
    \pi(x) &\defeq 0.62 \omega^2 \theta - 0.61\theta^3 + 10.14\theta, \nonumber \\
    \gamma(r) &\defeq 0.8r, \quad r \in \mathbb{R}_{\geqslant 0}. \nonumber 
\end{align} 
It is noteworthy that the synthesized DTCBF $h$ and its friend $\pi$ are independent of $x_c$ and $v_c$. The synthesized DTCBF is depicted in Figure \ref{fig:sec6:example-cart-pole}. To illustrate the level of conservatism in the proposed method for this numerical case study, we also plot the maximal polyhedral controlled invariant set for the system \eqref{eq:sec6:linear-system} with $\mathbb{U}$ using the MPT3 toolbox \cite{MPT3}, without considering \mbox{the safe set $\mathcal{S}$} since it cannot be defined within this toolbox. Specifically, the MPT3 toolbox computes the maximal polyhedral controlled invariant sets for linear systems considering polyhedral state and input constraints. By comparing the maximal polyhedral controlled invariant set with the zero-superlevel set of the synthesized DTCBF in Figure \ref{fig:sec6:example-cart-pole}, it can be observed that our proposed method synthesizes a DTCBF whose zero-superlevel set is nearly maximal for the system \eqref{eq:sec6:linear-system} with $\mathbb{U}$ and the safe set $\mathcal{S}$.

\subsection{Nonlinear System}
\begin{figure}[ht]\centering 
\includegraphics[width=0.42\textwidth]{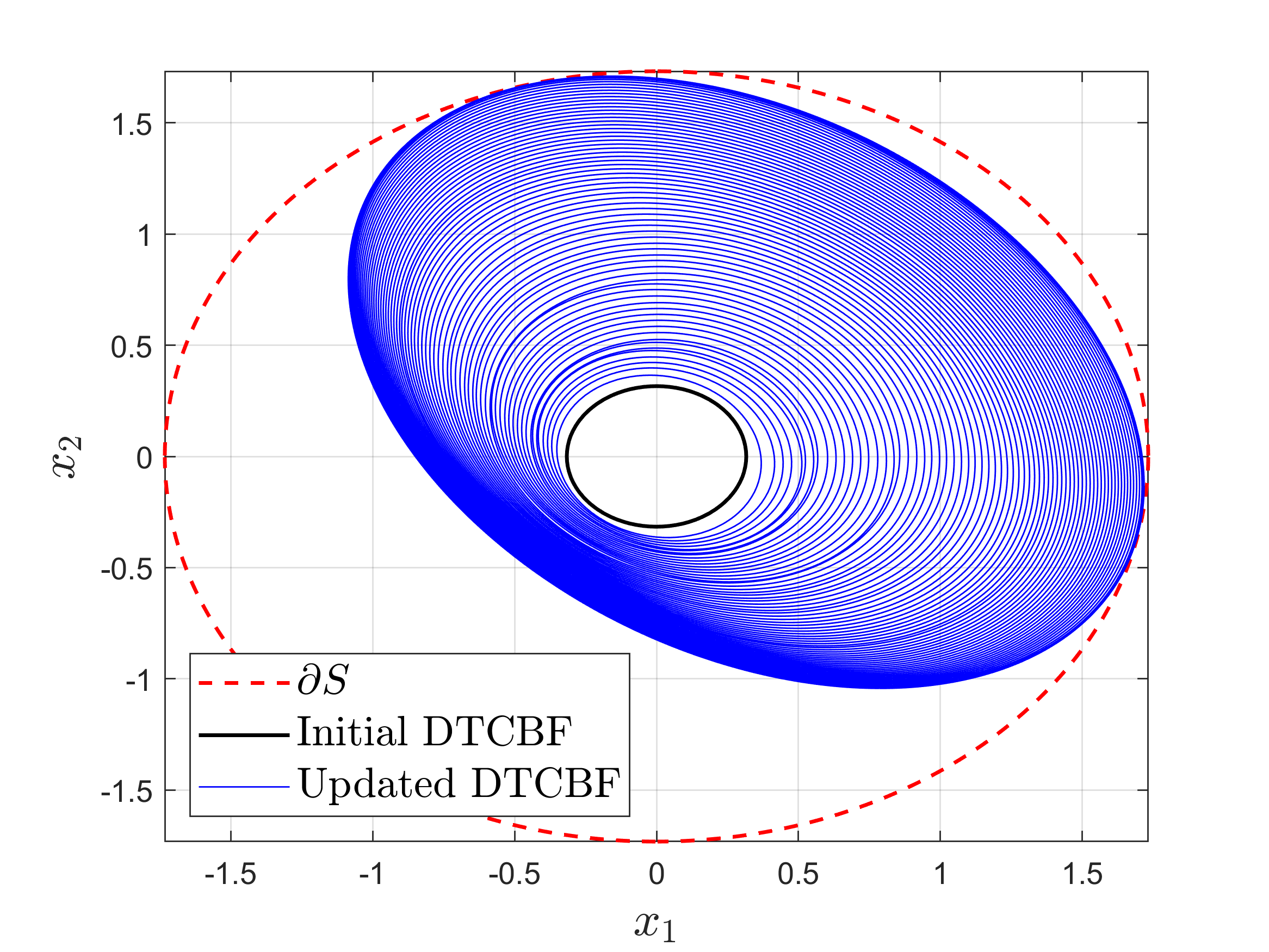}
    \caption{The proposed synthesis method is applied to the nonlinear system \eqref{eq:sec6:nonlinear-system} with the control admissible set $\mathbb{U}$ and the safe set \eqref{eq:sec6:safe-set} (the interior and boundary of the red circle) to iteratively synthesize a DTCBF-triple, starting from the initial DTCBF \eqref{eq:sec6:initial-guess} (black circle). Each of these blue curves represents the updated DTCBF after each iteration.}
    \label{fig:sec6:example-nonlinear}
\end{figure}

We discretize the continuous-time system discussed in \cite{Wang2023a} using the forward Euler method and with sampling \mbox{time $T_s = \unit[1]{s}$}, leading to 
\begin{align}
    \begin{bmatrix}
        x_1^+ \\
        x_2^+
    \end{bmatrix} &= \begin{bmatrix}
        x_1 + x_2T_s \\
        x_2 + \bigl(x_1 + \frac{1}{3}x_1^3 + x_2\bigr)T_s
    \end{bmatrix} + \nonumber \\
    &~ \begin{bmatrix}
        \bigl(x_1^2 + x_2 + 1\bigr)T_s & 0 \\
        0 & \bigl(x_2^2 + x_1 + 1\bigr)T_s
    \end{bmatrix} \begin{bmatrix}
        u_1 \\
        u_2
    \end{bmatrix}, \label{eq:sec6:nonlinear-system}
\end{align}
where $u \defeq [u_1 ~u_2]^{\T} \in \mathbb{U}$ with
\begin{align}
   \mathbb{U} \defeq \{ u \in \mathbb{R}^2 \mid -1.5 \leqslant u_1 \leqslant 1.5, ~ -1.5 \leqslant u_2 \leqslant 1.5 \}.
    \nonumber
\end{align}
In this example, the safe set $\mathcal{S}$ is defined as the disk
\begin{align} \label{eq:sec6:safe-set}
    \mathcal{S} \defeq \{x \in \mathbb{R}^2 \mid x_1^2 + x_2^2 - 3 \leqslant 0 \}.
\end{align}
We aim to synthesize a DTCBF-triple $(h, \gamma, \pi)$ according to the objectives of Problem \ref{prob:sec3:synthesis} for the system \eqref{eq:sec6:nonlinear-system} with the control admissible set $\mathbb{U}$ and the safe set \eqref{eq:sec6:safe-set} by applying the proposed method. In the initialization step, \mbox{$H \defeq \Pi_1 \hspace{-0.07cm} \defeq \hspace{-0.03cm} \Pi_2 \hspace{-0.07cm} \defeq \hspace{-0.03cm} [1 \hspace{0.14cm} x_1 \hspace{0.14cm}x_2 \hspace{0.14cm}x_1x_2 \hspace{0.14cm}x_1^2 \hspace{0.14cm}x_2^2]^{\T}$} are considered as vectors of monomials for the DTCBF $h$ and the control policy $\pi$, parameterized as in \eqref{eq:sec4:parameterized-DTCBF-general} and \eqref{eq:sec4:parameterized-control-policy}. Regarding $\gamma \in \mathcal{K}^{\leqslant \mathrm{id}}_{\infty}$ parameterized as in \eqref{eq:sec4:parameterized-gamma}, assume that we aim to maximize the value of $\gamma_0 \in (0,1] \subset \mathbb{R}$. Additionally, a small circle centered at the
origin is considered as an initial DTCBF, given by 
\begin{align} \label{eq:sec6:initial-guess}
    h^{(0)}(x) \defeq - x_1^2 - x_2^2 + 0.1.
\end{align} 
By applying the synthesis method, after around 60 iterations, with each iteration taking an average time of $\unit[5]{s}$ (total time about $\unit[5]{min}$), a DTCBF-triple $(h, \gamma, \pi)$ is \mbox{synthesized} iteratively as 
\begin{align}
    h(x) &\defeq -0.183x_1^2 - 0.124x_1x_2 - 0.189x_2^2 +  0.156x_1 + \nonumber \\
    &\qquad  0.164x_2 + 0.269, \nonumber \\
    \pi_1(x) &\defeq 0.139x_1^2 + 0.312x_1x_2 + 0.103x_2^2 - 0.681x_1 - \nonumber \\ 
    &\qquad   0.686x_2 + 0.211, \nonumber \\
    \pi_2(x) &\defeq 0.035x_1^2 + 0.324x_1x_2 + 0.159x_2^2 - 0.702x_1 - \nonumber \\
    &\qquad  0.877x_2 + 0.208, \nonumber \\
    \gamma(r) &\defeq r, \quad r \in \mathbb{R}_{\geqslant 0}. \nonumber
\end{align}
As shown in Figure \ref{fig:sec6:example-nonlinear}, the zero-superlevel set of the synthesized DTCBF cannot be further enlarged due to the restrictiveness of the safe set $\mathcal{S}$. 


\section{Conclusion and Future Work} \label{sec:sec8}
In this paper, we have proposed an alternating-descent approach based on Sum-of-Squares (SOS) programming to synthesize polynomial discrete-time Control Barrier Functions (DTCBFs) and corresponding polynomial control policies for discrete-time control-affine polynomial systems with input constraints and semi-algebraic safe sets. We have applied our approach to numerical case studies illustrating its efficacy. In our future work, we aim to synthesize \textit{robust} DTCBFs to guarantee safety for discrete-time systems subject to disturbances. 


\bibliographystyle{plain}        
\bibliography{autosam}           

\begin{thebibliography}{10}

\bibitem{Agrawal2017a}
A.~Agrawal and K.~Sreenath.
\newblock Discrete control barrier functions for safety-critical control of discrete systems with application to bipedal robot navigation.
\newblock In {\em Robotics: Science and Systems}, pages 1--10, 2017.

\bibitem{Ames2019a}
A.~D. Ames, S.~Coogan, M.~Egerstedt, G.~Notomista, K.~Sreenath, and P.~Tabuada.
\newblock Control barrier functions: Theory and applications.
\newblock In {\em 18th European Control Conference (ECC)}, pages 3420--3431, 2019.

\bibitem{Ames2014a}
A.~D. Ames, J.~W. Grizzle, and P.~Tabuada.
\newblock Control barrier function based quadratic programs with application to adaptive cruise control.
\newblock In {\em 53rd IEEE Conference on Decision and Control (CDC)}, pages 6271--6278, 2014.

\bibitem{chen2021}
Y.~Chen, M.~Jankovic, M.~Santillo, and A.~D. Ames.
\newblock Backup control barrier functions: Formulation and comparative study.
\newblock In {\em 60th IEEE Conference on Decision and Control (CDC)}, pages 6835--6841, 2021.

\bibitem{clark2022}
A.~Clark.
\newblock A semi-algebraic framework for verification and synthesis of control barrier functions, 2022.
\newblock \textit{arXiv preprint arXiv:2209.00081}.

\bibitem{LCP-book}
R.~W. Cottle, J.-S. Pang, and R.~E. Stone.
\newblock {\em The Linear Complementarity Problem}.
\newblock Philadelphia: Society for Industrial and Applied Mathematics, 2009.

\bibitem{dai2022}
H.~Dai and F.~Permenter.
\newblock {Convex synthesis and verification of control-Lyapunov and barrier functions with input constraints}.
\newblock In {\em American Control Conference (ACC)}, pages 4116--4123, 2023.

\bibitem{Giannis-TAC}
G.~Delimpaltadakis, J.~Cortés, and W.P.M.H. Heemels.
\newblock Continuous approximations of projected dynamical systems via control barrier functions.
\newblock {\em IEEE Transactions on Automatic Control}, 70(1):681--688, 2025.

\bibitem{Freire2023a}
V.~Freire and M.~M. Nicotra.
\newblock {Systematic Design of Discrete-Time Control Barrier Functions Using Maximal Output Admissible Sets}.
\newblock {\em IEEE Control Systems Letters}, 7:1891--1896, 2023.

\bibitem{RDTCBF-literature}
V.~Freire and M.~M. Nicotra.
\newblock Building robust control barrier functions from robust maximal output admissible sets.
\newblock In {\em 63rd IEEE Conference on Decision and Control (CDC)}, pages 8171--8177, 2024.

\bibitem{Gurriet2020}
T.~Gurriet, M.~Mote, A.~Singletary, P.~Nilsson, E.~Feron, and A.~D. Ames.
\newblock A scalable safety critical control framework for nonlinear systems.
\newblock {\em IEEE Access}, 8:187249--187275, 2020.

\bibitem{MPT3}
M.~Herceg, M.~Kvasnica, C.N. Jones, and M.~Morari.
\newblock Multi-parametric toolbox 3.0.
\newblock In {\em European Control Conference (ECC)}, pages 502--510, 2013.

\bibitem{Zamani-2020}
P.~Jagtap, S.~Soudjani, and M.~Zamani.
\newblock Formal synthesis of stochastic systems via control barrier certificates.
\newblock {\em IEEE Transactions on Automatic Control}, 66(7):3097--3110, 2021.

\bibitem{Kang2023}
S.~Kang, Y.~Chen, H.~Yang, and M.~Pavone.
\newblock Verification and synthesis of robust control barrier functions: Multilevel polynomial optimization and semidefinite relaxation.
\newblock In {\em 62nd IEEE Conference on Decision and Control (CDC)}, pages 8215--8222, 2023.

\bibitem{katriniok2023}
A.~Katriniok, E.~Shakhesi, and W.P.M.H. Heemels.
\newblock {Discrete-time control barrier functions for guaranteed recursive feasibility in nonlinear MPC: An application to lane merging}.
\newblock In {\em 62nd IEEE Conference on Decision and Control (CDC)}, pages 3776--3783, 2023.

\bibitem{lindemann2024}
L.~Lindemann, A.~Robey, L.~Jiang, S.~Das, S.~Tu, and N.~Matni.
\newblock Learning robust output control barrier functions from safe expert demonstrations.
\newblock {\em IEEE Open Journal of Control Systems}, 3:158--172, 2024.

\bibitem{sostools}
A.~Papachristodoulou, J.~Anderson, G.~Valmorbida, S.~Prajna, P.~Seiler, P.~A. Parrilo, M.~M. Peet, and D.~Jagt.
\newblock {\em {SOSTOOLS}: Sum of squares optimization toolbox for {MATLAB}}, 2021.
\newblock \textit{arXiv preprint arXiv:1310.4716}.

\bibitem{SOS2005}
A.~Papachristodoulou and S.~Prajna.
\newblock A tutorial on sum of squares techniques for systems analysis.
\newblock In {\em Proceedings of the American Control Conference}, pages 2686--2700 vol. 4, 2005.

\bibitem{S-procedure}
M.~Putinar.
\newblock Positive polynomials on compact semi-algebraic sets.
\newblock {\em Indiana University Mathematics Journal}, 42(3):969--984, 1993.

\bibitem{Hilbert1888}
B.~Reznick.
\newblock Some concrete aspects of \text{Hilbert}'s 17th problem.
\newblock {\em Contemporary Mathematics, American Mathematical Society}, 253:251--272, 2000.

\bibitem{Shakhesi2025}
E.~Shakhesi, W.P.M.H. Heemels, and A.~Katriniok.
\newblock Synthesis of quadratic discrete-time control barrier functions for polynomial systems based on sum-of-squares programming.
\newblock Accepted in 13th IFAC Symposium on Nonlinear Control Systems. \url{https://heemels.tue.nl/assets/content/papers/ShaHee_NOLCOS25a.pdf}.

\bibitem{shakhesi2024b}
E.~Shakhesi, W.P.M.H. Heemels, and A.~Katriniok.
\newblock Optimization-based verification of discrete-time control barrier functions: A branch-and-bound approach.
\newblock In {\em 63rd IEEE Conference on Decision and Control (CDC)}, pages 3632--3637, 2024.

\bibitem{Wang2023a}
H.~Wang, K.~Margellos, and A.~Papachristodoulou.
\newblock {\em Assessing Safety for Control Systems Using Sum-of-Squares Programming}, pages 207--234.
\newblock Springer Nature Switzerland, Cham, 2023.

\bibitem{Wang2023b}
H.~Wang, K.~Margellos, and A.~Papachristodoulou.
\newblock Safety verification and controller synthesis for systems with input constraints.
\newblock {\em IFAC-PapersOnLine}, 56(2):1698--1703, 2023.
\newblock 22nd IFAC World Congress.

\bibitem{wang2024convexcodesign}
H.~Wang, K.~Margellos, A.~Papachristodoulou, and C.~De Persis.
\newblock Convex co-design of control barrier function and safe feedback controller under input constraints, 2024.
\newblock \textit{arXiv preprint arXiv:2403.11763}.

\bibitem{wei2023}
T.~Wei, S.~Kang, W.~Zhao, and C.~Liu.
\newblock Persistently feasible robust safe control by safety index synthesis and convex semi-infinite programming.
\newblock {\em IEEE Control Systems Letters}, 7:1213--1218, 2023.

\bibitem{wu2023neural}
J.~Wu, A.~Clark, Y.~Kantaros, and Y.~Vorobeychik.
\newblock Neural {Lyapunov} control for discrete-time systems.
\newblock In {\em Advances in Neural Information Processing Systems}, volume~36, pages 2939--2955. Curran Associates, Inc., 2023.

\bibitem{Zhang2024}
Y.~Yang, Y.~Zhang, W.~Zou, J.~Chen, Y.~Yin, and S.~E.~Li.
\newblock Synthesizing control barrier functions with feasible region iteration for safe reinforcement learning.
\newblock {\em IEEE Transactions on Automatic Control}, 69(4):2713--2720, 2024.

\bibitem{Zeng2021b}
J.~Zeng, Z.~Li, and K.~Sreenath.
\newblock Enhancing feasibility and safety of nonlinear model predictive control with discrete-time control barrier functions.
\newblock In {\em 60th IEEE Conference on Decision and Control (CDC)}, pages 6137--6144, 2021.

\bibitem{Zeng2021a}
J.~Zeng, B.~Zhang, and K.~Sreenath.
\newblock Safety-critical model predictive control with discrete-time control barrier function.
\newblock In {\em American Control Conference (ACC)}, pages 3882--3889, 2021.

\bibitem{Zhang2023}
H.~Zhang, Z.~Li, H.~Dai, and A.~Clark.
\newblock Efficient sum of squares-based verification and construction of control barrier functions by sampling on algebraic varieties.
\newblock In {\em 62nd IEEE Conference on Decision and Control (CDC)}, pages 5384--5391, 2023.

\end{thebibliography}



\appendix
\section{Proof of Proposition \ref{prop:sec5:prop1}} \label{sec:appendix-1}
By Lemma \ref{lemma:sec5:det-nonnegative}, \eqref{eq:prop:bi-A-1} holds, if and only if for all $x \in \mathbb{R}^n$, 
\begin{align}
    \mathrm{Y}_{i,1}(x) - \pi_i^2(x) &\geqslant 0, \label{eq:appendix-A-proof1} \\
    \mathrm{Y}_{i,1}(x) &\geqslant 0, \label{eq:appendix-A-proof2}
\end{align}
$i \in \{1,2, \hdots, m\}$. 
Based on the generalized S-procedure (Lemma \ref{lemma:sec4:S-procedure}), \eqref{eq:appendix-A-proof1} ensures that for all \mbox{$x \in \mathcal{C}^{(k-1)}$} such that $a^{(k-1)}_{i,i}(x) \leqslant 0$, 
\begin{align} \label{eq:appendix-A-proof3}
    \widetilde{\pi}_{i,i}(x) - \pi_i^2(x) \geqslant 0.
\end{align}
Note that satisfying \eqref{eq:appendix-A-proof1} implies that \eqref{eq:appendix-A-proof2} also holds, so the condition \eqref{eq:appendix-A-proof2} does not introduce conservatism.
Similarly, \eqref{eq:prop:bi-A-2} implies that for all \mbox{$x \in \mathcal{C}^{(k-1)}$} such that $a^{(k-1)}_{i,i}(x) \geqslant 0$,
    it holds that $\widetilde{\pi}_{i,i}(x) \leqslant 0$, and thus 
\begin{align} \label{eq:appendix-A-proof4}
    \pi_i^2(x) - \widetilde{\pi}_{i,i}(x) \geqslant 0.
\end{align}
In conclusion, the satisfaction of \eqref{eq:appendix-A-proof3} and \eqref{eq:appendix-A-proof4} ensures that \eqref{eq:prop1-final} holds for all $x \in \mathcal{C}^{(k-1)}$. \qed
\section{Proof of Proposition \ref{prop:sec5:prop2}} \label{sec:appendix-2}
Similar to the proof of Proposition \ref{prop:sec5:prop1}, satisfying \eqref{eq:prop:bi-B-1} and \eqref{eq:prop:bi-B-2} implies that for all $x \in \mathcal{C}^{(k-1)}$ such that \mbox{$a^{(k-1)}_{i,j}(x) \leqslant 0$}, $i \in \{1,2, \hdots, m-1\}$, $j \in \{i+1, i+2, \hdots, m\}$, it holds that 
\begin{align}
    \Theta_{i,j,1}(x) &\geqslant \pi^2_i(x), \label{eq:appendix-B-proof1} \\
    \Theta_{i,j,2}(x) &\geqslant \pi^2_j(x). \label{eq:appendix-B-proof2}
\end{align}
Moreover, following the generalized S-procedure (Lemma \ref{lemma:sec4:S-procedure}), the SOS constraints \eqref{eq:prop:bi-B-3} and \eqref{eq:prop:bi-B-4} imply that for all $x \in \mathcal{C}^{(k-1)}$ such that $a^{(k-1)}_{i,j}(x) \leqslant 0$,
\begin{align}
    \Theta_{i,j,3}(x) &\geqslant 0 \geqslant -(\pi_i(x) - \pi_j(x))^2, \label{eq:appendix-B-proof3} \\
    2\widetilde{\pi}_{i,j}(x) &\geqslant \Theta_{i,j,1}(x) + \Theta_{i,j,2}(x) + \Theta_{i,j,3}(x). \label{eq:appendix-B-proof4}
\end{align}
By summing up \eqref{eq:appendix-B-proof1}, \eqref{eq:appendix-B-proof2}, and \eqref{eq:appendix-B-proof3}, it follows that for all $x \in \mathcal{C}^{(k-1)}$ such that $a^{(k-1)}_{i,j}(x) \leqslant 0$,
\begin{align} \label{eq:appendix-B-proof5}
    \Theta_{i,j,1}(x) + \Theta_{i,j,2}(x) + \Theta_{i,j,3}(x) \geqslant 2\pi_i(x)\pi_j(x).
\end{align}
In conclusion, satisfying \eqref{eq:appendix-B-proof4} and \eqref{eq:appendix-B-proof5} implies that for all \mbox{$x \in \mathcal{C}^{(k-1)}$} such that $a^{(k-1)}_{i,j}(x) \leqslant 0$, it holds that 
\begin{align}
    \widetilde{\pi}_{i,j}(x) - \pi_i(x)\pi_j(x) \geqslant 0.
\end{align} \qed

\end{document}